\documentclass[11pt]{ip-journal}
\usepackage{pdfsync}
\usepackage{fullpage}
\usepackage{amsmath}
\usepackage{amssymb}
\usepackage{verbatim}
\usepackage[linktocpage]{hyperref}

\newtheorem{df}{Definition}[section]
\newtheorem{thm}[df]{Theorem}

\newtheorem{rem}[df]{Remark}

\newtheorem{lem}[df]{Lemma}

\newtheorem{cor}[df]{Corollary}
\newtheorem{conj}[df]{Conjecture}

\newcommand{\eq}{eqnarray*}

\newcommand{\lig}{\mathfrak{g}}
\newcommand{\lik}{\mathfrak{k}}

\title{Minimal surfaces for Hitchin representations}

\author{Song Dai\textsuperscript{1}}

\address{Song Dai\\
Center for Applied Mathematics at Tianjin University\\
Tianjin University\\
No.92 Weijinlu Nankai District\\
Tianjin\\
P.R.China 300072\\}

\email{song.dai@tju.edu.cn}

\author{Qiongling Li\textsuperscript{2}}

\address{Qiongling Li\\
Centre for Quantum Geometry of Moduli Spaces (QGM)\\
Department of Mathematics, Aarhus University\\
Ny Munkegade 118, Bldg. 1530\\
8000 Aarhus C\\
Denmark\\}

\address{Department of Mathematics\\
California Institute of Technology\\
1200 East California Boulevard\\
 Pasadena, CA 91125\\}

\email{qiongling.li@gmail.com}

\date{}

\begin{document}

\begin{abstract}
Given a reductive representation $\rho: \pi_1(S)\rightarrow G$, there exists a $\rho$-equivariant harmonic map $f$ from the universal cover of a fixed Riemann surface $\Sigma$ to the symmetric space $G/K$ associated to $G$. If the Hopf differential of $f$ vanishes, the harmonic map is then minimal. In this paper, we investigate the properties of immersed minimal surfaces inside symmetric space associated to a subloci of Hitchin component: the $q_n$ and $q_{n-1}$ cases. First, we show that the pullback metric of the minimal surface dominates a constant multiple of the hyperbolic metric in the same conformal class and has a strong rigidity property. Secondly, we show that the immersed minimal surface is never tangential to any flat inside the symmetric space. As a direct corollary, the pullback metric of the minimal surface is always strictly negatively curved. In the end, we find a fully decoupled system to approximate the coupled Hitchin system.
\end{abstract}
\maketitle

\footnotetext[1]{The author is supported by NSFC grant No. 11601369.}

\footnotetext[2]{Corresponding author, supported in part by the center of excellence grant `Center for Quantum Geometry of Moduli Spaces' from the Danish National Research Foundation (DNRF95).}

\section{Introduction}
For a closed, connected, oriented surface $S$ of genus $g\geq 2$ and a reductive Lie group $G$, consider the representation variety $Rep(\pi_1(S),G)=Hom(\pi_1(S),G)//G$.  For $G=PSL(2,\mathbb{R})$, there are two connected components of $Rep(\pi_1(S),PSL(2,\mathbb{R}))$ are identified with Teichm\"uller space; and the representations in these components are called Fuchsian. For a general real split Lie group, using the unique irreducible representation $PSL(2,\mathbb{R})\rightarrow G$, we can single out a component of $Rep(\pi_1(S),G)$, the connected component containing representations factors through Fuchsian representations, called Hitchin component for $G$. In particular, we denote the Hitchin component for $PSL(n,\mathbb{R})$ as Hit$_n$.

Fix a Riemann surface structure $\Sigma$ on $S$. By the work of Donaldson \cite{Donaldson} and Corlette \cite{Corlette}, given a reductive representation $\rho: \pi_1(S)\rightarrow G$, there exists a $\rho$-equivariant harmonic map $f$ from the universal cover of $\Sigma$ to the symmetric space $G/K$ associated to $G$, where $K$ is the maximal compact subgroup of $G$. If the representation $\rho$ is a Hitchin representation, Sanders \cite{HitchinImmersion} showed that the corresponding equivariant harmonic map is an immersion. In particular, if the Hopf differential of the harmonic map $f$ vanishes, then the harmonic map is conformal and hence is a minimal immersion.

Understanding such equivariant minimal immersions $f: \tilde{\Sigma}\rightarrow G/K$ for Hitchin representations is the main goal of this paper. 

To achieve this goal, we need to use the tool of Higgs bundles. By the results of Hitchin \cite{Hitchin87} and Simpson \cite{Simpson88}, given a polystable $G$-Higgs bundle $(E,\phi)$, there exists a unique Hermitian metric $h$ (compatible with $G$-structure) satisfying Hitchin equation
\begin{\eq}
F^{\nabla^{h}}+[\phi,\phi^{*_{h}}]=0
\end{\eq}
giving rise to a flat connection $D=\nabla^{h}+\phi+\phi^{*_{h}}$, the corresponding holonomy $\rho:\pi_1(S)\rightarrow G$, and a $\rho$-equivariant harmonic map $f: \tilde{\Sigma}\rightarrow G/K$.

Hitchin \cite{Hitchin92} gives an explicit description for the Hitchin component for $SL(n,\mathbb{R})$ in terms of Higgs bundles. Explicitly, the bundle of the Higgs bundles is $E=K^{\frac{n-1}{2}}\oplus K^{\frac{n-3}{2}}\oplus\cdots\oplus K^{\frac{1-n}{2}}$ and the Higgs field is explicitly parametrized by the holomorphic differentials $(q_2,q_3,\cdots,q_n)\in\bigoplus\limits_{i=2}^n H^0(K^i)$. In particular, $q_2$ is the Hopf differential of the corresponding harmonic map.

Replacing the quadratic differential by varying the Riemann surface choice, Labourie \cite{LabourieEnergy} considered the Hitchin map from the total space of vector bundle over Teichm\"uller space with fiber at $\Sigma$ as $(0,q_3,\cdots,q_n)$ to $\text{Hit}_n$ and showed that this map is surjective. Therefore, understanding all minimal surfaces arising from the Hitchin representations will eventually give properties of the representations, not depending on the choice of Riemann surface. In the same paper, Labourie conjectured that the Hitchin map is also injective.  When $n=3,$ Labourie \cite{LabourieCubic} and Loftin \cite{LoftinAffine} independently proved the conjecture using affine geometry. Using different methods, Labourie \cite{LabourieCyclic} then proved the conjecture for Hitchin representations into all rank 2 real split Lie groups: $PSL(3,\mathbb{R})$, $PSp(4,\mathbb{R}), G_2,$ and $PSL(2,\mathbb{R})\times PSL(2,\mathbb{R})$. Equivalently, given any Hitchin representation $\rho$ into a rank 2 real split Lie group $G$, there is a unique minimal surface inside $G/K$ that is $\rho$-equivariant.

We restrict to consider the Higgs fields parametrized by $(0,\cdots,0,q_n)$ and $(0,\cdots,q_{n-1},0)$:
\begin{\eq}
q_n ~ \text{case:}
\left(
\begin{array}{cccccc}
0 & & & & & q_n\\
1 & 0 & & & &\\
& 1 & 0 & & &\\
& & & \ddots & &\\
& & & & 0 &\\
& & & & 1 & 0
\end{array}
\right),&\qquad&
q_{n-1} ~ \text{case:}
\left(
\begin{array}{cccccc}
0 & & & & q_{n-1} & 0\\
1 & 0 & & & & q_{n-1}\\
& 1 & 0 & & &\\
& & & \ddots & &\\
& & & & 0 &\\
& & & & 1 & 0
\end{array}
\right).
\end{\eq}
In particular, $q_n=0$ or $q_{n-1}=0$ gives the base Fuchsian point. From the surjectivity of the Hitchin map, for any Hitchin representation, we may find a complex structure on $S$, such that $q_2=0$ in the corresponding Higgs bundle. Hence, if we vary the choice of Riemann surface, these two families of Higgs bundles give the whole Hitchin component for $PSL(3,\mathbb{R})$, $PSp(4,\mathbb{R})$, and $G_2$.

If $n=2$, the Higgs bundles of the $q_2$ case parametrize the whole Teichm\"uller component. If $n=3$, the $q_{n-1}$ case gives the embedding of Teichm\"uller component inside Hitchin component for $PSL(3,\mathbb{R})$.  But in thses two casese, the Hopf differential of the harmonic map does not vanish and hence the harmonic map is not conformal. However, we remark that in these two cases, almost all of our results can be applied. For simplicity of language, from now on, whenever referring to the $q_n, q_{n-1}$ cases, we don't include the case that $q_n$ (or $q_{n-1}$) is $q_2$.

Such Higgs bundles of the $q_n$ case were initially studied by Baraglia \cite{Bar}, which he called cyclic Higgs bundle. Later Collier \cite{BrianThesis} considered other Higgs bundles under finite order automorphisms, in particular $q_{n-1}$ case (also see \cite{CL14}). These families of Higgs bundles possess particular nice properties. In both cases, the Hermitian metric $h$ solving Hitchin equations is diagonal, i.e., $h=(h_1,h_2,\cdots, h_2^{-1},h_1^{-1})$.

We are ready to explain our main results for the $q_n,q_{n-1}$ cases. In general, these are some families in Hitchin component. But from the above argument, we remark that our results in fact hold for the Hitchin component for $PSL(3,\mathbb{R})$, $PSp(4,\mathbb{R})$, and $G_2$.
\\

$\bullet$  \textbf{Metric Domination.} We show that the pullback metric $g_f$ of the minimal immersion $f$ dominates the base Fuchsian metric $g_{\text{Fuchsian}}$.
\begin{thm}\label{dom1}
For the $q_n$, $q_{n-1}$ cases
\begin{\eq}
g_f=2n tr(\phi\phi^*)\geq g_{\text{Fuchsian}}=\frac{1}{6}(n^4-n^2)g_0.
\end{\eq}
Moreover, if equality holds at one point, then $q_n=0$, or $q_{n-1}=0$ respectively, which implies it is base Fuchsian.
\end{thm}

For the precise definitions above involving the $q_n$, $q_{n-1}$ cases, see Section \ref{problem}.

If we integrate the pullback metric, this is closely related to the Morse function considered in Hitchin \cite{Hitchin92}, which plays an important role to determine the topology of the representation variety for $PSL(n,\mathbb{R})$. The Morse function on moduli space of polystable Higgs bundles is defined as: \begin{\eq} f(E,\phi)=\int_{
\Sigma}\text{tr}(\phi\phi^*)  \sqrt{-1}dz\wedge d\bar z.\end{\eq} Hitchin \cite{Hitchin92} showed that in the Hitchin component, the only minimum is the base Fuchsian point. That is, consider any Higgs bundle $(E,\phi)$ in Hitchin component parametrized by $(q_2,q_3,\cdots,q_n)$, then\begin{\eq}
 f(E,\phi)\geq \text{base Fuchsian case}=\text{topological quantity}.
\end{\eq}
Equality holds if and only if it is base Fuchsian.

When $q_2=0$, the harmonic map is a minimal immersion and the Morse function is in fact the area of the minimal surface up to a constant. Then \begin{\eq}
\text{MinArea}({\rho})\geq\text{Area}(\rho_{\text{Fuchsian}})=\text{topological quantity},
\end{\eq}
where $\text{MinArea}({\rho})$ is the minimum of the area of the pullback metric going through all the $\rho$-equivariant immersion. Equality holds if and only if it is the base Fuchsian representation $\rho_{\text{Fuchsian}}$. Labourie \cite{LabourieCyclic} pointed out that this is also a corollary of the entropy result for Hitchin representations by Potrie-Sambarino \cite{HitchinEntropy}. For the Hitchin component for rank 2 Lie group, this is reproved by Labourie \cite{LabourieCyclic} which he called area rigidity formula.

In fact, this area rigidity formula by Labourie inspired us to be interested in this metric domination question. We conjecture that this domination property also holds for all minimal surfaces arising from Hitchin representations.
\begin{conj} (Metric Domination Conjecture)\label{DominationConjecture}
Consider Higgs bundle $(E,\phi)$ in Hitchin component parametrized by $(0,q_3,\cdots,q_n)$. On the surface $\Sigma$, the pullback metric $g_f$ of the minimal immersion $f$ satisfies
\begin{\eq}
g_f\geq g_{\text{Fuchsian}}=\frac{1}{6}(n^4-n^2)g_0.
\end{\eq}
If the equality holds at one point, then it holds at every point.
\end{conj}

\begin{rem} For $q_2\neq 0$, in other words, the harmonic map $f$ is not conformal. Instead of the pullback metric $g_f$, we consider the $(1,1)$ part of $g_f$, $g_{f}^{(1,1)}=2n\text{tr}(\phi\phi^{*_h})$ which is conformal to the base Fuchsian metric. From our proof of Theorem \ref{dom1}, the above domination theorem is also true in non-minimal surface case for the lower rank case. More precisely, for Hitchin representations in $PSL(2,\mathbb{R})$ parametrized by $q_2$, we have
 \begin{\eq}
g_{f}^{(1,1)}\geq g_{\text{Fuchsian}}=2g_0.
\end{\eq}
Moreover, if equality holds at one point, then $q_2=0$, that is, the base Fuchsian case. \end{rem}

The domination question above is comparing Hitchin representations with the base Fuchsian points. If instead, we consider the other non-Hitchin components, it won't make sense to compare them with Fuchsian points. At first, the harmonic map is not necessarily an immersion, hence the pullback metric can never dominate a hyperbolic metric pointwise. A more desired question is to ask whether the representations dominate some particular special representations which are similar to the status of Fuchsian points in Hitchin components.

It is interesting to compare the domination question here with another type of domination question in terms of length spectrum. By work of Deroin and Tholozan \cite{DominationFuchsian} and Gu\'{e}ritaud, Kassel and Wolff \cite{Kassel}, for representations of $\pi_1(S)$ into $PSL(2,\mathbb{R})$,  any Fuschian representation dominates some non-Fuchsian representation and any non-Fuchsian representation can be dominated by some Fuchsian one. In a similar spirit, Lee and Zhang \cite{CollarLemma} conjectured for any Hitchin representation $\rho:\pi_1(S)\rightarrow PSL(n, R),$ there is a Fuchsian representation whose length spectrum strictly dominated by $\rho$. Tholozan \cite{TwoEntropy} proved this conjecture for $n=3$ case and mentioned that Labourie pointed out the conjecture of Lee and Zhang cannot hold anymore for $n\geq 4$. The contradiction comes from Hitchin representations in $PSp(2k, R)$ and $PSO(k, k + 1)$.  Therefore, in the same paper, Tholozan made some modification of the conjecture by changing Fuchsian representation by representations into $SO(n,n+1)$ or $Sp(2n,\mathbb{R})$. Our Conjecture \ref{DominationConjecture} here is in fact weaker than the conjecture by Lee and Zhang since the length spectrum of the pullback back metric is larger than the length spectrum of the representation which is Lipschitz to the distance inside the symmetric space. So even though our conjecture is already true for Hitchin representations into rank 2 Lie groups: $PSL(3,\mathbb{R}), PSp(4,\mathbb{R})$, and $G_2$, it does not imply the conjecture of Lee and Zhang.
\\

$\bullet$ \textbf{Negative Curvature.} We describe how the immersed minimal surface sits inside the symmetric space by showing that the minimal immersion is never tangential to any flat inside the symmetric space $G/K$.
\begin{thm}\label{neg}
For the $q_n$, $q_{n-1}$ cases, we have the following results. \\
(1) The Hitchin equation never decouples: for every point on $\Sigma$, \begin{\eq} F^{\nabla^{h}}\neq 0, \quad [\phi,\phi^{*_h}]\neq 0.\end{\eq}
(2) The sectional curvature $K_{G/K}(\sigma)$ is strictly negative, where $\sigma$ is the tangent space of the image of $f$. Geometrically, the minimal immersion is never tangential to any flat inside the symmetric space.\\
(3) On each line bundle $K^{\frac{n+1-2k}{2}}$, $1\leq k\leq n$, the Chern form of $h_k$: \begin{\eq}
\sqrt{-1}\Theta_{h_k}=\sqrt{-1}\overline{\partial}\partial\log h_k \end{\eq} is strictly positive if $n+1-2k>0$; zero if $n+1-2k=0$; strictly negative if $n+1-2k<0$.
\end{thm}

\begin{rem}
The phenomenon in Part (1) is in contrast to the asymptotic behaviour of Hitchin equation proved in Collier and Li \cite {CL14}: along the ray $tq_n$ (or $tq_{n-1}$), the Hitchin equation decouples as $t\rightarrow\infty$: \begin{\eq} F^{\nabla^{h}}\rightarrow 0, \quad [\phi,\phi^{*_h}]\rightarrow 0.\end{\eq}  This asymptotic behavior is generalized by Mochizuki \cite{Takuro} to a much more general family of Higgs bundles.
\end{rem}
From Gauss equation, it is easy to see that the curvature of the pullback metric $g_f$ is always non-positive. Then moreover, as a corollary of Theorem \ref{neg}, we obtain the curvature is strictly negative.
\begin{cor}\label{Negative} For the $q_n$, $q_{n-1}$ cases, the sectional curvature of the immersed minimal surface is strictly negative.\end{cor}
We conjecture that this phenomenon is true for all minimal immersions arising from Hitchin representations.
\begin{conj}(Negative Curvature Conjecture) For the Hitchin representation parametrized by $(0,q_3,\cdots,q_n)$, the minimal immersion is never tangential to any flat inside the symmetric space. And as a corollary, the sectional curvature of immersed minimal surface is strictly negatively curved.
\end{conj}
We apply our result to estimate the entropy of Hitchin representations. Given $\rho$ a Hitchin representation into $G$ and select a point in the symmetric space $p\in G/K$. The volume entropy of $\rho$ is defined as \begin{\eq}
h(\rho):=\underset{R\rightarrow \infty}{\text{lim sup }}\frac{\text{log}(\#|\{\gamma\in \pi_1(S)|d(p,\rho(\gamma)(p))\leq R\}|)}{R},
\end{\eq} where $d$ is the distance in $G/K$. Lots of progress have been made on the volume entropy of Hitchin representations. Potrie and Sambarino \cite{HitchinEntropy} showed that for any Hitchin representation $\rho$, one has $h(\rho)\leq 1$ and the equality holds only if $\rho$ is Fuchsian. Zhang \cite{Zhang3,ZhangN} constructed certain sequences of Hitchin representations
along which $h(\rho)\rightarrow 0.$ For $n=3$, Nie \cite{XinNie} showed that the entropy of Hitchin representations parametrized by $tq_3$ goes to zero as $t\rightarrow\infty$.
Sanders \cite{HitchinImmersion} showed that for a Hitchin representation, the curvature of the pullback metric $g_f$ satisfies\begin{\eq} \frac{1}{\text{Area}  (g_f)}\int_{\Sigma}\sqrt{-K_{g_f}}d V_{g_f}\leq h(\rho).\end{\eq} And he used this inequality to show $h(\rho)>0$. The entropy of Hitchin representation then satisfies:
\begin{\eq}
h(\rho)\geq \underset{\Sigma}{\min}\{\sqrt{-K_{g_f}}\}
\end{\eq}
By Corollary \ref{Negative}, the sectional curvature of immersed minimal surface $K_{g_f}<0$ for the $q_n,q_{n-1}$ cases, we obtain
\begin{cor}$\underset{\Sigma}{\min}\{\sqrt{-K_{g_f}}\}$ provides a positive lower bound for the volume entropy $h(\rho)$ of Hitchin representations for the $q_n, q_{n-1}$ cases.\end{cor}

In fact, since $K_{g_f}\leq K_{G/K}$, we can also use $\underset{\Sigma}{\min}\{\sqrt{-K_{G/K}}\}$ as a weaker lower bound. The term $\underset{\Sigma}{\min}\{\sqrt{-K_{G/K}}\}$ is very interesting and involves more analytical terms. We hope to show a more quantitative estimate on this term in future work. \\

$\bullet$ \textbf{Coupled Hitchin System vs Decoupled System.} We investigate the coupled Hitchin system in the $q_n,q_{n-1}$ cases. Coupled equations are generally hard to study. If there is some way to decouple the system, it will be substantially easier to solve the system. In this paper, we find a fully decoupled system such that the solutions of the decoupled system approximate the solutions of coupled Hitchin system of equations. The decoupled system are formed of single scalar equations as follows. For $n$ even, set $n=2m$; for $n$ odd, set $n=2m+1$.
\begin{thm}\label{scalar1}
There exists a unique Hermitian metric $u_k, v_k$ respectively on the holomorphic line bundle $K^{\frac{n+1-2k}{2}}$ satisfying, locally
\begin{eqnarray*}
\triangle \log u_k+u_k^{-\frac{2}{n+1-2k}}-(u_k^2|q_{n}|^2)^{\frac{1}{2k-1}}&=&0,\quad 1\leq k\leq m, \\
\triangle \log v_k+v_k^{-\frac{2}{n+1-2k}}-(v_k^2|2q_{n-1}|^2)^{\frac{1}{2k-2}}&=&0,\quad 2\leq k\leq m.
\end{eqnarray*}
And the following estimates hold,
\begin{eqnarray*}
\max\{|q_n|^{\frac{2\alpha_k}{n}},(\alpha_kg_0)^{\alpha_k}\}\leq &u_k^{-1}&<((\underset{\Sigma}{\max}|q_n|_{g_0}^{\frac{2}{n}}+\alpha_k)g_0)^{\alpha_k}, \quad 1\leq k\leq m,\\
\max\{|2q_{n-1}|^{\frac{2\alpha_k}{n-1}}, (\alpha_kg_0)^{\alpha_k}\}\leq &v_k^{-1}&<((\underset{\Sigma}{\max}|2q_{n-1}|_{g_0}^{\frac{2}{n-1}}+\alpha_k)g_0)^{\alpha_k}, \quad 2\leq k\leq m,
\end{eqnarray*}
where $\alpha_k=\frac{n+1-2k}{2}$ and $g_0$ is the Hermitian hyperbolic metric on $K^{-1}$.
\end{thm}
\begin{rem}
Such equations are natural generalization of familiar vortex-like equations. In particular, in the case $n=2,k=1$, this is the harmonic equation from surface to surface with given Hopf differential $q_2$. In the case $n=3, k=1$, this gives Wang's equation \cite{Wang} and the solution is the Blaschke metric for the hyperbolic affine sphere. Dumas and Wolf \cite{DumasWolf} solved such equations for the case when $n$ is general and $k=1$ on the complex plane.
\end{rem}

\begin{rem}
Surprisingly, there is a geometric interpretation for the above system of equations. Let $\sigma_k=u_k^{-\frac{2}{n+1-2k}}$ be an Hermitian metric on the Riemann surface $\Sigma$. The equation for $\sigma_k$ is
\begin{\eq} \frac{n+1-2k}{2}K_{\sigma_k}=-1+|q_{n}|_{\sigma_k}^{\frac{2}{2k-1}},
\end{\eq}
where $K_{\sigma_k}$ is the curvature of $\sigma_k$. Roughly speaking, the curvature of $\sigma_k$ differs from the curvature of hyperbolic metric by the norm of $q_n$.
\end{rem}

We show that $u_k$ (or $v_k$) is an upper approximate of $h_k$. Let $u_k^t$ (or $v_k^t$) be the solution for $tq_n$ (or $tq_{n-1}$) in Theorem \ref{scalar1}.
\begin{thm}\label{comparison1} For $m\geq 2$, suppose $q_n$ (or $q_{n-1}$) is not zero. Then for each compact set $K\subset \Sigma$ away from zeros of $q_n$ (or $q_{n-1}$), there is a positive constant $C=C(K)$ independent of $t$ such that
\begin{\eq}
\text{$q_{n}$ case:}& \quad u_k^t(1-C||tq_n||^{-1})\leq h_k^t<u_k^t, \quad 1\leq k\leq m,\\
\text{$q_{n-1}$ case:}& \quad v_k^t(1-C||tq_{n-1}||^{-1})\leq h_k^t< v_k^t, \quad 2\leq k\leq m.
\end{\eq}
\end{thm}

\begin{rem}
To prove the left direction of the above theorem, we make use the asymptotic behavior of $h_k$ established by Collier and Li \cite{CL14}: for example, for the $q_n$ case, 
\begin{\eq}
h_k=|q_n|^{-\frac{n+1-2k}{n}}(1+O(||q_n||^{-1}))
\end{\eq}
away from the zeros of $q_n$. And our results in fact improve the estimate of Collier and Li from the other direction by the estimate of $u_k,v_k$ in Theorem \ref{scalar1}.
\end{rem}

\textbf{Structure of the article.} The article is organized as follows. In Section \ref{pre}, we recall some fundamental results about Higgs bundles. In particular, we recall the Donaldson-Uhlenbeck-Yau correspondence and the explicit relation between Higgs bundles and related harmonic maps for further calculation. We fix some notations at the end of Section \ref{pre}. In Section \ref{problem}, we set up the main object and describe the $q_n$ case and $q_{n-1}$ case. In Section \ref{domi}, we show the metric domination Theorem \ref{dom1}. In Section \ref{shape}, we describe the shape of minimal surface inside the symmetric space by proving Theorem \ref{neg}. In Section \ref{decouple}, we find a decoupled system to bound the Hitchin coupled system and  prove Theorem \ref{scalar1} and Theorem \ref{comparison1}.

\subsection*{Acknowledgement} The authors wish to thank Mike Wolf, Andy Sanders, Ian McIntosh, Daniele Alessandrini and Brian Collier for helpful discussions and comments. The second author acknowledges support from U.S. National Science Foundation grants DMS 1107452, 1107263, 1107367 ``RNMS: GEometric structures And Representation varieties" (the GEAR Network).

\section{Preliminaries and notations}\label{pre}
In this section, we recall some results in (principal) Higgs bundle and harmonic map. A good reference is the thesis of Baraglia \cite{Bar}, Section 2.1.

Let $\Sigma$ be a closed Riemann surface with genus $g\geq 2$. Let $\pi_1=\pi_1(\Sigma,p)$ be the fundamental group. Let $G$ be a complex semisimple Lie group, $K$ be a compact real form of $G$. Let $\mathfrak{g}$, $\mathfrak{k}$ be the corresponding Lie algebra. Then $\mathfrak{g}=\mathfrak{k}\otimes \mathbb{C}$ as a real Lie algebra. Let $\theta$ be the anti-linear involution of $\mathfrak{g}$ induced from the conjugation of $\mathbb{C}$. Let $B_{\mathfrak{g}}$ be the complex Killing form of $\mathfrak{g}$ and $B_{\mathfrak{k}}$ be the real Killing form of $\mathfrak{k}$. Then $B_{\mathfrak{g}}$ is the complex linear extension of $B_{\mathfrak{k}}$. So without confusing, we just use $B$ to denote Killing form of $\lig$ or $\lik$. Then $H(X,Y):=-B(X,\theta{Y})$ gives an Hermitian metric on $\lig$. Then we have an orthogonal decomposition $\lig=\lik\oplus\lik^{\bot}$ with respect to $H$, where $\lik^{\bot}=\sqrt{-1}\lik$. This decomposition also gives a direct sum as $Ad_{K}$ module. And notice that $[\lik^{\bot},\lik^{\bot}]\subseteq \lik$.

For example, $G=SL(n,\mathbb{C})$, $K=SU(n)$. Then $B(X,Y)=2n\text{tr}_{\mathbb{C}}(XY)$, where $X,Y\in sl(n,\mathbb{C})$.

We can establish a correspondence between the following two moduli spaces.

\textbf{Moduli Space A:} Equivalent classes $[\rho]$ of reductive representations, $\rho:\pi_1\to G$. Here reductive means that the induced representation on $\lig$ is a direct sum of irreducible representations. For two representations $\rho_1,\rho_2$, they are equivalent if and only if there exists $g\in G$, such that $\rho_1=g\rho_2 g^{-1}$.

\textbf{Moduli Space B:} Gauge equivalent classes of $G$-Higgs bundles satisfying Hitchin equations. More precisely, let $P$ be a principal $G$ bundle over $\Sigma$. Suppose we have a $K$ principal bundle reduction $i: P_K\hookrightarrow P$, or equivalently, a section of $P\times_{l}G/K$, where $l$ is the left multiplication. Denote $AdP^{c}=P\times_{Ad_{G}}\lig$. This reduction also gives an identification
\begin{\eq}
AdP^{c}=P\times_{Ad_{G}}(\lik\otimes\mathbb{C})=P_K\times_{Ad_{K}}(\lik\otimes\mathbb{C})=(P_K\times_{Ad_{K}}\lik)\otimes\mathbb{C}.
\end{\eq}
Notice that $\theta$ and $B$ are both $Ad_{K}$ invariant on $\lik\otimes\mathbb{C}$. So as a complex vector bundle $AdP^{c}$, we can define $\theta$ and $B$ on $AdP^{c}$, and then Hermitian metric $H$. Let $X\in AdP^{c}$, denote $X^{*}=-\theta(X)$ for the sake that the adjoint of $ad_X$ is $ad_{-\theta(X)}$ with respect to $H$. Let $\phi$ be a section of $AdP^{c}\otimes K$, where $K$ is the canonical line bundle (don't confuse with the compact subgroup $K$). Locally, suppose $\phi=Xdz$, then we define $\phi^{*}=X^*d\bar{z}\in AdP^{c}\otimes \bar{K}$ and $[\phi,\phi^*]=[X,X^*]dz\wedge d\bar{z}$. Let $A$ be a principal connection on $P_K$, $\nabla_A$ be the corresponding connection on associated bundles of $P_K$. Denote $\bar{\partial}_{A}=\nabla_{A}^{0,1}\otimes 1+1\otimes\bar{\partial}$, which gives a holomorphic structure on $AdP^{c}\otimes K$. Denote $F_A$ be the curvature of $A$. Then the Hitchin equations are given by
\begin{eqnarray}
F_A+[\phi,\phi^*]&=&0,\label{h1}\\
\bar{\partial}_{A}\phi&=&0\label{h2}.
\end{eqnarray}
For a $G$-Higgs bundle satisfying Hitchin equations, we just mean the data $(\Sigma,G,K,P,P_K,\phi,A)$, or briefly $(\phi, A)$, satisfying equations (\ref{h1}) and (\ref{h2}). Notice that $A+\phi+\phi^*$ gives a principal connection on $P$, and equation (\ref{h1}) is equivalent to this connection being flat. Given $\Sigma, G, K$, for two $G$-Higgs bundles satisfying Hitchin equations, $(\phi_1, A_1)$ and $(\phi_2, A_2)$, they are equivalent if and only if there exists a $K$-gauge transformation $\alpha$ over $\Sigma$ (automatically $G$-gauge transformation), such that $\alpha^*\phi_2=\phi_1, \alpha^*A_2=A_1$.

Before we discuss the relation between Moduli Space A and Moduli Space B, we first consider the relation between (vector) Higgs bundle and $G$-Higgs bundle satisfying Hitchin equations.
\begin{df}
Let $\Sigma$ be a Riemann surface, $E$ be a holomorphic vector bundle over $\Sigma$. Let $\phi$ be a holomorphic section of $End(E)\otimes K$. We call $(E,\phi)$ a Higgs bundle over $\Sigma$.
\end{df}

Given a $G$-Higgs bundle satisfying the Hitchin equations and suppose $G$ acts on $\mathbb{C}^{n}$. Then the associated bundle gives a complex vector bundle $E$, and $\nabla_{A}^{0,1}$ gives a holomorphic structure of $E$. By definition, $\phi$ is a section of $End(E)\otimes K$. And by equation $(\ref{h2})$, $\phi$ is holomorphic.

Conversely, under some assumptions, one can obtain a $G$-Higgs bundle satisfying Hitchin equations from Higgs bundle. Denote $\mu_{E}=\frac{\deg(E)}{\text{rank}(E)}$ be the slope of $E$, where $\deg(E)$ is the degree of $E$, and $\text{rank}(E)$ is the complex rank of $E$. We call $(E,\phi)$ is stable if for any proper $\phi$-invariant holomorphic subbundle $F$, $\mu_{F}<\mu_{E}$. We call $(E,\Phi)$ is polystable if $(E,\Phi)$ is a direct sum of stable Higgs bundles. Let $G=SL(n,\mathbb{C})$, $K=SU(n)$. We have the following result.
\begin{thm}(Hitchin \cite{Hitchin87} and Simpson \cite{Simpson88})\label{dh}
Let $(E,\phi)$ be a polystable Higgs bundle with structure group $SL(n,\mathbb{C})$. Then there exists a unique Hermitian metric $h$ on $E$ compatible with the $SL(n,\mathbb{C})$ structure, such that
\begin{eqnarray}
F^{\nabla^h}+[\phi,\phi^{*_h}]=0, \label{curvature equation}
\end{eqnarray}
where $F^{\nabla^h}$ is the Chern connection of $h$, in local holomorphic trivialization,
\begin{\eq}
F^{\nabla^h}=\overline{\partial}(h^{-1}\partial h),
\end{\eq}
and $\phi^{*_{h}}$ is the adjoint of $\phi$ with respect to $h$, in the sense that
\begin{\eq}
h(\phi(u),v)=h(u,\phi^{*_{h}}(v)),\quad u,v\in \Gamma(E)
\end{\eq}
in local frame, $\phi^{*_{h}}=\overline{h}^{-1}\overline{\phi}^{\top}\overline{h}$, and the bracket is the commutator of $End(E)$. We regard both $F^{\nabla_h}$ and $[\phi,\phi^{*_h}]$ as sections of $(K\wedge\overline{K})\otimes End(E)$. This gives rise to a flat connection $\nabla^{h}+\phi+\phi^{*_{h}}$.
\end{thm}
Notice that the Hermitian metric $h$ gives a reduction to $SU(n)$ bundle, and ${\nabla^h}$ gives a principal $SU(n)$ connection. Then clearly we obtain a $G$-Higgs bundle satisfying Hitchin equations.

Now we establish the correspondence between Moduli Space A and Moduli Space B. Given a $G$-Higgs bundle satisfying Hitchin equations, since $A+\phi+\phi^*$ is a flat principal connection on $P$, the monodromy gives a representation $\rho:\pi_1\to G$. One can show that $\rho$ is reductive. And the monodromy descends to the equivalent classes.

Given a reductive representation $\rho$, we have a associated principal $G$ bundle $\tilde{\Sigma}\times_{\rho}G$, denoted as $P$, where $\tilde{\Sigma}$ is the universal cover of $\Sigma$, regarded as the $\pi_1$ principal bundle over $\Sigma$. Now we want to find a reduction from $P$ to a principal $K$ bundle $P_K$. The reduction should satisfy some condition which will be clarified below. Let $i: P_K\hookrightarrow P$ be the reduction, which is equivalent to a $\rho$-equivariant map $f: \tilde{\Sigma}\to G/K$. Notice that the Maurer-Cartan form $\omega$ of $G$ gives a flat connection on $P$, we still use $\omega$ to denote the connection. Consider $i^{*}\omega$, which is a $\lig$ value one form on $P_K$. Decomposing $i^{*}\omega=A+\Phi$ from $\lig=\lik\oplus\lik^{\bot}$, where $A$ is $\lik$ valued and $\Phi$ is $\lik^{\bot}$ valued. Then $A$ is a principal connection on $P_K$ and $\Phi$ is a section of $T^*M\otimes (P_{K}\times_{Ad_{K}}\lik^{\bot})$. Consider
\begin{\eq}
T^*M\otimes (P_{K}\times_{Ad_{K}}\lik^{\bot})\otimes\mathbb{C}=(T^*M\otimes\mathbb{C})\otimes (P_{K}\times_{Ad_{K}}\lik^{\bot}\otimes\mathbb{C})=(K\oplus\bar{K})\otimes AdP^{c}.
\end{\eq}
Regard $\Phi$ as a section of $T^*M\otimes (P_{K}\times_{Ad_{K}}\lik^{\bot})\otimes\mathbb{C}$, and consider the decomposition $(K\oplus\bar{K})\otimes AdP^{c}$. Let $\Phi=\phi+\phi^*$, where $\phi$ is a section of $K\otimes AdP^{c}$ and $\phi^*$ is a section of $\bar{K}\otimes AdP^{c}$. Notice that $\theta(\Phi)=-\Phi$, so $\phi^*=-\theta(\phi)$. Now we obtain all the desired data. Equation (\ref{h1}) follows from the flatness of $\omega$. By direct calculation (see \cite{Bar}), equation (\ref{h2}) is equivalent to $f$ being harmonic, where the conformal metric on $\tilde{\Sigma}$ is induced from the complex structure and the metric on $G/K$ is induced from the Killing form of $\lig$. So our requirement on $f$ is just that $f$ is $\rho$-equivariant and harmonic. From Donaldson \cite{Donaldson} and Corlette \cite{Corlette}, the existence of such $f$ is equivalent to $\rho$ being reductive. And $f$ is unique up to the composition of the centraliser of $\rho(\pi_1)$. We see the above construction descends to the equivalence classes. So we have built the bijection between Moduli Space A and Moduli Space B.\\

So far, given a reductive representation $\rho:\pi_1\to G$, we have a $\rho$-equivariant harmonic map $f:\tilde{\Sigma}\to G/K$. Let $\hat{g}$ be the metric on $G/K$ induced from Killing form of $\lig$. More precisely, let $B$ be the Killing form of $\lig$. Consider the tangent bundle $T(G/K)=G\times_{Ad_K}\lik^{\bot}$ over $G/K$. Since $B$ is $Ad_K$-invariant and positive definite on $\lik^{\bot}$, we have a well defined Riemannian metric on $G/K$. Locally, suppose $U$ is a neighborhood of $p\in G/K$, $X,Y\in T_pU$. Let $\tilde{U},\tilde{p},\tilde{X},\tilde{Y}$ be a lift to $G$. Using the Maurer-Cartan form, we have a decomposition of $T_{\tilde{p}}G$ according to $\lig=\lik\oplus\lik^{\bot}$. Then $\hat{g}(X,Y)=B(\tilde{X}^{\bot},\tilde{Y}^{\bot})$, where $\tilde{X}^{\bot}$ is the $\lik^{\bot}$ component of $\tilde{X}$.\\

The Higgs bundle and the corresponding harmonic map are related as follows.

Locally, choose a lift from $\Sigma$ to $\tilde{\Sigma}$ and a lift from $G/K$ to $G$. And lift $f$ to $\tilde{f}$ as a map from $\tilde{\Sigma}$ to $G$. For $x\in\Sigma$, we see that, \begin{\eq}(\tilde{f}_*\tilde{X})^{\bot}=\Phi(X),\end{\eq} where $X\in T_{x}\Sigma$. Hence,
\begin{\eq}
(\tilde{f}_*\widetilde{\frac{\partial}{\partial z}})^{\bot}=\phi(\frac{\partial}{\partial z}),\quad(\tilde{f}_*\widetilde{\frac{\partial}{\partial\bar{z}}})^{\bot}=\phi^*(\frac{\partial}{\partial \bar{z}}).
\end{\eq}
We consider the pullback metric $g_f$ on $\Sigma$, $g_f=\pi_{*}f^{*}\hat{g}$, where $\pi$ is the covering map $\pi: \tilde{\Sigma}\to \Sigma$. Since $f$ is $\rho$-equivariant and $\hat{g}$ is $G$-invariant, $g_f$ is well defined. Then $\forall X,Y\in T\Sigma$,
\begin{\eq}
g_f(X,Y)=B(\Phi(X),\Phi(Y)).
\end{\eq}
 In particular, for $G=SL(n,\mathbb{C})$, $K=SU(n)$, we have \begin{\eq}\text{Hopf}(f)=g_f^{2,0}=2n\text{tr}_{\mathbb{C}}(\phi\phi).\end{\eq}If $\text{Hopf}(f)=0$, then as a section of $K\otimes\bar{K}$, the Hermitian metric is \begin{\eq}g_f=g_f^{1,1}=2n\text{tr}_{\mathbb{C}}(\phi\phi^{*_{h}}).\end{\eq} (In fact, $\phi\phi^{*_{h}}$ is real.) And the corresponding Riemannian metric is $g_f+\bar{g}_f$.
\\

Now we fix some notations used throughout this paper.

Let $g$ be an Hermitian metric on $K^{-1}$, where $K$ is the canonical line bundle. We can also regard $g$ as a section of $K\otimes\bar{K}$. In local coordinate, $g=gdz\otimes d\bar{z}$. We abuse the same notation $g$ to denote both the metric and the local function. Similarly, let $h$ be an Hermitian metric on $K^{-l}$, denoted by $h=hdz^{\otimes l}\otimes d\bar{z}^{\otimes l}$. Notice that as a local function, $h$ is not globally well defined, but if we set $h=ag^{l}$, then $a$ is a globally defined function.

 Let $q_n$ be a $n$-differential. Locally, denote $q_n=q_ndz^ n$. Define a local function $|q_n|^2=q_n\bar{q}_n$, which is corresponding to the Hermitian metric on $K^{-n}$. Here are some notations used later.\\
(1) $|q_n|^2_g$ as the square of the norm of $q_n$ with respect to the metric $g$, which is globally defined, in local coordinate, $|q_n|^2_g=|q_n|^2g^{-n}$;\\
(2) $||q_n||=\int_{\Sigma} |q_n|^{\frac{2}{n}}$, which will be used in Section \ref{decouple};\\
(3) $g_0$ as the unique Hermitian hyperbolic metric compatible with the complex structure;\\
(3) $g_{\text{Fuchsian}}$ as the Hermitian metric corresponding to the base Fuchsian point, a multiple of $g_0$;\\
(4) $\triangle=\partial_{z}\partial_{\overline{z}}=\frac{1}{4}(\frac{\partial^{2}}{\partial x^{2}}+\frac{\partial^{2}}{\partial y^{2}})$, which is locally defined;\\
(5) $\triangle_{g}=g^{-1}\triangle$, noting that $\triangle_{g}$ is globally defined.

\section{Rewriting Hitchin equations in two subclasses for Hitchin representations}\label{problem}
In Section \ref{pre}, we relate the moduli space of reductive representations $\rho\in Hom(\pi_1,SL(n,\mathbb{C}))$ to the moduli space of $SL(n,\mathbb{C})$ Higgs bundles satisfying Hitchin equations. Under this bijection, we can also describe the moduli space of reductive representations $\rho\in Hom(\pi_1,SL(n,\mathbb{R}))$ in the setting of Higgs bundles. In fact, in \cite{Hitchin92}, Hitchin gives a parametrization of the Hitchin component of the $SL(n,\mathbb{R})$-Higgs bundle moduli space using Higgs bundles of the form $(E,\phi)$ as follows. Let
\begin{\eq}
E=S^{n-1}(K^{\frac{1}{2}}\oplus K^{-\frac{1}{2}})=K^{\frac{n-1}{2}}\oplus K^{\frac{n-3}{2}}\oplus\dots\oplus K^{-\frac{n-3}{2}}\oplus K^{-\frac{n-1}{2}}
\end{\eq}
be the $(n-1)$'st symmetric power and the Higgs field $\phi$ is explicitly parametrized by $(q_2,q_3,\dots, q_n)\in\bigoplus\limits_{j=2}^nH^0(\Sigma, K^j).$ The embedded copy of Teichm\"uller space comes from setting $q_3=\dots=q_n=0$.

We restrict to consider two subclasses of Higgs bundles in Hitchin component: the $q_n$, $q_{n-1}$ cases.

\textbf{$q_n$ Case:}
The Higgs field $\phi$ is a holomorphic section of $End(E)\otimes K $ of the form
\begin{displaymath}
\left(
\begin{array}{cccccc}
0 & & & & & q_n\\
1 & 0 & & & &\\
& 1 & 0 & & &\\
& & & \ddots & &\\
& & & & 0 &\\
& & & & 1 & 0
\end{array}
\right),
\end{displaymath}
where $q_n$ is a holomorphic section of $K^n$. From \cite{Bar}, the Hermitian metric $h$ is given by
\begin{\eq}
h=\text{diag}(h_1,h_2,\dots,h^{-1}_{2},h^{-1}_{1}),
\end{\eq}
where $h_{k}$ is an Hermitian metric on $K^{\frac{n-2k+1}{2}}$, i.e., $h_{k}$ is a positive definite smooth section of $K^{-\frac{n-2k+1}{2}}\otimes \overline{K}^{-\frac{n-2k+1}{2}}$. We also denote $K\otimes \overline{K}$ as $|K|^{2}$.

Sanders \cite{HitchinImmersion} showed that the corresponding equivariant harmonic map $f$ for Hitchin representations is an immersion. By direct calculation, the Hopf differential $\text{Hopf}(f)=2n\text{tr}_{\mathbb{C}}(\phi^2)$ vanishes and hence $f$ is a minimal immersion.
Equation (\ref{curvature equation}) gives the following system of equations.

For $n=2m$ even, 
\begin{eqnarray*}
\triangle \log h_1+h_1^{-1}h_2-h_1^2|q_{n}|^2&=&0,\nonumber\\
\triangle \log h_k+h_k^{-1}h_{k+1}-h_{k-1}^{-1}h_{k}&=&0,\quad 2\leq k\leq m-1\\
\triangle \log h_m+h_m^{-2}-h_{m-1}^{-1}h_{m}&=&0,\nonumber
\end{eqnarray*}
where $\triangle=\partial_{z}\partial_{\overline{z}}=\frac{1}{4}(\frac{\partial^{2}}{\partial x^{2}}+\frac{\partial^{2}}{\partial y^{2}})$ is the coordinate Laplacian.\\
If $q_n=0$, by the uniqueness of the solution, one may solve the system of equations as
\begin{\eq}
h^{-1}_{k}h_{k+1}=\frac{1}{2}k(n-k)g_{0},\quad 1\leq k\leq m-1,\qquad h^{-2}_{m}=\frac{1}{8}n^2g_0,
\end{\eq}
where $g_0$ is the hyperbolic Hermitian metric on $K^{-1}$.
\\
The pullback metric $g_f$ of the minimal immersion on the surface is
\begin{\eq}
g_f=2n\text{tr}_{\mathbb{C}}(\phi\phi^{*_{h}})=2n(h_1^2|q_{n}|^2+2\sum_{k=1}^{m-1}h_{k}^{-1}h_{k+1}+h_{m}^{-2}).
\end{\eq}

For $n=2m+1$ odd, the situation is similar. Notice that there is a trivial bundle $K^0$ in the middle of $E$, we see that the metric on $K^0$ is $1$, i.e., the standard Hermitian metric on $\mathbb{C}$. Then the last equation of the system of equations becomes
\begin{\eq}
\triangle \log h_{m}+h_{m}^{-1}-h_{m-1}^{-1}h_{m}&=&0.\nonumber
\end{\eq}
If $q_n=0$, then
\begin{\eq}
h^{-1}_{k}h_{k+1}=\frac{1}{2}k(n-k)g_{0},\quad 1\leq k\leq m-1,\qquad h^{-1}_m=\frac{1}{8}(n^2-1)g_0.
\end{\eq}
The last term in the pullback metric $g_f$ of the minimal immersion is $2h_{m}^{-1}$ instead of $h_{m}^{-2}$.

\textbf{$q_{n-1}$ Case:}
The Higgs field $\phi$ is a holomorphic section of $End(E)\otimes K$ of the form
\begin{displaymath}
\left(
\begin{array}{cccccc}
0 & & & & q_{n-1} & 0\\
1 & 0 & & & & q_{n-1}\\
& 1 & 0 & & &\\
& & & \ddots & &\\
& & & & 0 &\\
& & & & 1 & 0
\end{array}
\right),
\end{displaymath}
where $q_{n-1}$ is a holomorphic section of $K^{n-1}$. From \cite{CL14}, the Hermitian metric $h$ is given by
\begin{\eq}
h=\text{diag}(h_1,h_2,\dots,h^{-1}_{2},h^{-1}_{1}),
\end{\eq}
By direct calculation, the Hopf differential $\text{Hopf}(f)=2n\text{tr}_{\mathbb{C}}(\phi^2)$ vanishes and hence $f$ is a minimal immersion.
Equation (\ref{curvature equation}) gives the following system of equations.

For $n=2m$ even, 
\begin{\eq}
\triangle \log h_1+h_1^{-1}h_2-h_1 h_2|q_{n-1}|^2&=&0,\nonumber\\
\triangle \log h_2+h_2^{-1}h_3-h_1^{-1} h_2-h_1h_2|q_{n-1}|^2&=&0,\nonumber\\
\triangle \log h_k+h_k^{-1}h_{k+1}-h_{k-1}^{-1}h_{k}&=&0,\quad 3\leq k\leq m-1\\
\triangle \log h_m+h_m^{-2}-h_{m-1}^{-1}h_{m}&=&0.\nonumber
\end{\eq}
When $q_{n-1}=0$, $h_k$ is the same as in the $q_n$ case, $n=2m$.\\
The pullback metric of minimal immersion is
\begin{\eq}
g_f=2n\text{tr}_{\mathbb{C}}(\phi\phi^{*_{h}})=2n(2h_1h_2|q_{n-1}|^2+2\sum_{k=1}^{m-1}h_{k}^{-1}h_{k+1}+h_{m}^{-2}).
\end{\eq}
For $n=2m+1$ is odd, the last equation of the system of equations is
\begin{\eq}
\triangle \log h_m+h_m^{-1}-h_{m-1}^{-1}h_{m}&=&0.
\end{\eq}
When $q_{n-1}=0$, $h_k$ is the same as in the $q_n$ case, $n=2m+1$.
\\
The last term in the pullback metric $g_f$ is $2h_{m}^{-1}$ instead of $h_{m}^{-2}$.
\begin{rem}
For the situation $q_n=0$ in $q_n$ case (or $q_{n-1}=0$ in $q_{n-1}$ case), it is the base Fuchsian point. We denote the pullback metric $g_f$ in this situation as $g_{\text{Fuchsian}}$. It is computed in Section \ref{domi} that $g_{\text{Fuchsian}}=\frac{1}{6}(n^4-n^2)g_0$, where $g_0$ is the unique Hermitian hyperbolic metric.
\end{rem}

\section{Domination of pullback metric}\label{domi}
Let $g_0$ be the unique hyperbolic metric on Riemann surface $\Sigma$ with constant curvature $-1$ and compatible with the complex structure. Regarding $g_0$ as a section of $K\otimes\bar{K}$, we set $g_0=g_0dz\otimes d\bar{z}$. (If we regard $g_0$ as an Hermitian metric on $T\Sigma$, then $g_0=g_0(dz\otimes d\bar{z}+d\bar{z}\otimes dz)=2g_0(dx^2+dy^2)$.)

Locally, $g_{0}$ satisfies
\begin{\eq}
\triangle \log g_{0}-g_{0}=0,
\end{\eq}
where $\triangle=\partial_{z}\partial_{\overline{z}}=\frac{1}{4}(\frac{\partial^{2}}{\partial x^{2}}+\frac{\partial^{2}}{\partial y^{2}})$ is the coordinate Laplacian.

In this section, we show that the pullback metric of minimal immersion for the $q_n, q_{n-1}$ cases dominates the base Fuchsian metric $g_{\text{Fuchsian}}$. More precisely,
\begin{thm}\label{dom}
For the $q_n$, $q_{n-1}$ cases
\begin{\eq}
g_f\geq g_{\text{Fuchsian}}=\frac{1}{6}(n^4-n^2)g_0.
\end{\eq}
Moreover, if equality holds at one point, then $q_n=0$, or $q_{n-1}=0$ respectively, which implies $h_k$ is a suitable power of $g_0$, as described in Section \ref{problem}.
\end{thm}
\begin{rem}
In fact, we prove the rigidity result for every term in $g_f$. For example, for the $q_n$ case, $n=2m$, if one of the following equalities hold at one point,
\begin{\eq}
h^{-1}_{k}h_{k+1}\geq(km-\frac{1}{2}k^2)g_0,\quad 1\leq k\leq m-1,\qquad h^{-2}_m\geq\frac{1}{2}m^2g_0,
\end{\eq}
the rigidity result then holds.
\end{rem}
\textit{Proof of Theorem \ref{dom}.}
We only prove the theorem for the $q_n$ case, $n=2m$. For other cases, the proof is similar. Let $h_k=a_kg_0^{k-\frac{1}{2}-\frac{n}{2}},  1\leq k\leq m$. Let $|q_{n}|^{2}=a_0g_0^{n}$. Then $a_{k},  0\leq k\leq m$ are globally defined functions. They satisfy
\begin{\eq}
\triangle \log a_1+(\frac{1}{2}-\frac{n}{2})\triangle \log g_0+a_1^{-1}a_2g_0-a_{1}^{2}a_0g_0&=&0\\
\triangle \log a_k+(k-\frac{1}{2}-\frac{n}{2})\triangle \log g_0+(a_k^{-1}a_{k+1}-a_{k-1}^{-1}a_{k})g_0&=&0,\quad 2\leq k\leq m-1\\
\triangle \log a_m-\frac{1}{2}\triangle \log g_0+(a_m^{-2}-a_{m-1}^{-1}a_{m})g_0&=&0.
\end{\eq}
Then $g_f$ becomes $2n(a_1^2a_0+2\sum_{k=1}^{m-1}a_{k}^{-1}a_{k+1}+a_{m}^{-2})g_0.$
Let $b_{k}=a_{k}^{-1}a_{k+1}, 1\leq k\leq m-1$ and $b_{m}=a_{m}^{-2}$. Then
\begin{\eq}
a_k=b_{m}^{-\frac{1}{2}}b_{m-1}^{-1}\cdots b_{k}^{-1},\quad 1\leq k\leq m-1,\qquad a_m=b_m^{-\frac{1}{2}}.
\end{\eq}
So
\begin{\eq}
\triangle \log a_k&=&-\frac{1}{2}\triangle\log b_m-\triangle \log b_{m-1}-\cdots-\triangle\log b_k,\quad 1\leq k\leq m-1\\
\triangle \log a_m&=&-\frac{1}{2}\triangle\log b_m.
\end{\eq}
Plug these terms into the system of equations above, notice $\triangle\log g_0=g_0$, we have
\begin{\eq}
\frac{1}{2}\triangle \log b_m+\triangle \log b_{m-1}+\cdots+\triangle \log b_1-(\frac{1}{2}-m+b_1)g_0&\leq& 0\\
\frac{1}{2}\triangle \log b_m+\triangle \log b_{m-1}+\cdots+\triangle \log b_k-(k-\frac{1}{2}-m+b_k-b_{k-1})g_0&=& 0,\quad 2\leq k\leq m-1\\
\frac{1}{2}\triangle \log b_m-(-\frac{1}{2}+b_m-b_{m-1})g_0&=&0.
\end{\eq}
Then
\begin{eqnarray}
\triangle\log b_1+(1+b_2-2b_1)g_0&\leq& 0\nonumber\\
\triangle\log b_k+(1+b_{k+1}+b_{k-1}-2b_k)g_0&=& 0,\quad 2\leq k\leq m-1\label{b}\\
\triangle\log b_m+(1+2b_{m-1}-2b_m)g_0&=& 0.\nonumber
\end{eqnarray}
Let $x_{k}$ be a minimizer of $b_{k}$. Then
\begin{\eq}
2b_1&\geq& 1+b_{2}(x_1)\\
2b_{k}&\geq& 1+b_{k+1}(x_k)+b_{k-1}(x_k),\quad 2\leq k\leq m-1\\
2b_m&\geq& 1+2b_{m-1}(x_m).
\end{\eq}
First consider
\begin{\eq}
2b_{m-1}&\geq& 1+b_{m}(x_{m-1})+b_{m-2}(x_{m-1})\\
2b_m&\geq& 1+2b_{m-1}(x_m).
\end{\eq}
We see $b_m(x_{m-1})\geq 2+b_{m-2}(x_{m-1}).$
And then $2b_{m-1}\geq 3+2b_{m-2}(x_{m-1}).$
Then similarly, consider $k=m-1,m-2$. We obtain
\begin{\eq}
b_{m-1}(x_{m-2})&\geq& 4+b_{m-3}(x_{m-2})\\
2b_{m-2}&\geq& 5+2b_{m-3}(x_{m-2}).
\end{\eq}
Follow this procedure until we obtain
\begin{\eq}
b_{3}(x_{2})&\geq& 2(m-2)+b_{1}(x_{2})\\
2b_{2}&\geq& 2m-3+2b_{1}(x_2).
\end{\eq}
From $2b_1\geq 1+b_2(x_1)$, we have $b_2(x_1)\geq 2m-2.$ Then $2b_{1}\geq 2m-1.$ And then
\begin{\eq}
2b_k\geq (2m-1)+\cdots+(2m-(2k-1)), \quad 1\leq k\leq m.
\end{\eq}
So we obtain
\begin{\eq}
2b_k&\geq& 2km-k^{2},\quad 1\leq k\leq m.
\end{\eq}
Recall that $\sum_{k=1}^{m-1}k^{2}=\frac{1}{6}(m-1)m(2m-1). $ Hence finally,
\begin{\eq}
g_f=2n(h_1^2|q_{2m}|^2+2\sum_{k=1}^{m-1}h_{k}^{-1}h_{k+1}+h_{m}^{-2})\geq \frac{2n}{6}(4m^3-m)g_0=\frac{1}{6}(n^4-n^2)g_0.
\end{\eq}
For rigidity, suppose $b_{k}=km-\frac{1}{2}k^2$ holds at one point $p$ for some $k$, let $\log b_k=\Omega_k$, then from (\ref{b}),
\begin{\eq}
\triangle\Omega_k+2(km-\frac{1}{2}k^2-e^{\Omega_k})g_0\leq 0.
\end{\eq}
Then one may apply the strong maximum principle to finish the proof. More precisely, in a local coordinate chart,
\begin{\eq}
\triangle (\Omega_k-\log (km-\frac{1}{2}k^2))-2(\Omega_k-\log (km-\frac{1}{2}k^2))(\int_{0}^{1}e^{t\Omega_k+(1-t)\log (km-\frac{1}{2}k^2)}dt) g_0\leq 0.
\end{\eq}
This is from the fact
\begin{\eq}
e^{\Omega_k}-e^{\log (km-\frac{1}{2}k^2)}=\frac{d}{dt}\Big|_{0}^{1}e^{t\Omega_k+(1-t)\log (km-\frac{1}{2}k^2)}.
\end{\eq}
Notice that $\int_{0}^{1}e^{t\Omega_k+(1-t)\log (km-\frac{1}{2}k^2)}dt\geq 0$ and $\Omega_k-\log (km-\frac{1}{2}k^2)\geq 0$, then by the strong maximum principle \cite{MaximumPrincipal}, if $\Omega_k-\log (km-\frac{1}{2}k^2)=0$ at one point $p$, then $\Omega_k-\log (km-\frac{1}{2}k^2)= 0$ everywhere. Then $b_k$ equals $km-\frac{1}{2}k^2$ identically for all $1\leq k\leq m$. From the origin system of equations, we see that $q_{n}$ must be zero.
\qed

\section{Shape of minimal surface inside symmetric space}\label{shape}
In this section, we investigate the shape of minimal surface $\Sigma$ inside the symmetric space $G/K$. In particular, we show that the tangent space of $\Sigma$ is never tangential to any flat inside $G/K$. Recall $G=SL(n,\mathbb{C})$, $K=SU(n)$.
\begin{thm} \label{nevertangential} For the $q_n$, $q_{n-1}$ cases, we have the following results. \\
(1) The Hitchin equation never decouples: for every point on $\Sigma$, \begin{\eq} F^{\nabla^{h}}\neq 0, \quad [\phi,\phi^{*_h}]\neq 0.\end{\eq}
(2) The sectional curvature $K_{G/K}(\sigma)$ is strictly negative, where $\sigma$ is the tangent space of the image of $g_f$. Geometrically, the minimal immersion is never tangential to any flat inside the symmetric space.\\
(3) On each line bundle $K^{\frac{n+1-2k}{2}}$, $1\leq k\leq n$, the Chern form of $h_k$: \begin{\eq}
\sqrt{-1}\Theta_{h_k}=\sqrt{-1}\overline{\partial}\partial\log h_k\end{\eq} is strictly positive if $n+1-2k>0$; zero if $n+1-2k=0$; strictly negative if $n+1-2k<0$.
\end{thm}

Before we prove Theorem \ref{nevertangential}, an immediate corollary is as follows.
\begin{cor}\label{negative} For the $q_{n},q_{n-1}$ cases, the sectional curvature $K_{g_f}$ of $\Sigma$ equipped with the induced metric $g_f$ is strictly negative.
\end{cor}
\textit{Proof of Corollary \ref{negative}.}
Given the $\rho$-equivariant harmonic map $f:\widetilde{\Sigma}\rightarrow G/K$, we want to investigate that, as an immersed submanifold inside the symmetric space $G/K$, how $f(\widetilde{\Sigma})$ interacts with the symmetric space.
Let $\nabla$ denote the Levi-Civita connection of $G/K$ and $\nabla^T$ denote the component of $\nabla$ tangential to the image of $f$. Then the second fundamental form is the symmetric 2-tensor with values in the normal bundle given by
\begin{\eq}
II(X,Y)=\nabla_XY-\nabla^T_XY,
\end{\eq}
where $X,Y\in \Gamma(f^*T(G/K))$ are tangent to the image of $f$.  Let $\{e_1,e_2\}$ be an orthonormal basis such that
\begin{\eq}
II(e_1, e_2)=0.
\end{\eq}
The immersion is minimal if the trace of $II$ with respect to $g_f$ is vanishing. Therefore \begin{\eq}
II(e_1, e_1)+II(e_2,e_2)=0.
\end{\eq}

Recall the Gauss equation, let $X, Y$ be orthonormal vector fields on $\Sigma$, then
\begin{\eq}K_{g_f}(X, Y )=K_{G/K}(X, Y) + <II(X, X),II(Y, Y )>-|II(X, Y )|^2.\end{\eq}

Therefore the sectional curvature $K_{g_f}$ and $K_{G/K}$ measured in $\Sigma$ and in $G/K$ respectively are related by Gauss equation:
\begin{\eq}
K_{g_f}=K_{G/K}-\frac{1}{2}||II||^2_{g_f}.
\end{\eq}
By Theorem \ref{nevertangential}, $K_{G/K}<0.$ Therefore $K_{g_f}<0.$
\qed\\

\textit{Proof of Theorem \ref{nevertangential}.}
We first discover the following inequalities for the 2-tensors appearing in the Hitchin system of equations. We will prove Lemma \ref{inequalityfirst}, \ref{inequalitysecond} later.
\begin{lem}\label{inequalityfirst} For the $q_n$ case, on the whole surface $\Sigma$, if $n=2m$
\begin{\eq}
h_1^2|q_{n}|^2<h_1^{-1}h_2<h_2^{-1}h_3<\cdots<h_{k-1}^{-1}h_k<\dots<h_{m-1}^{-1}h_m<h_{m}^{-2}.
\end{\eq}
If $n=2m+1$, the last term is replaced by $h_{m}^{-1}.$
\end{lem}
\begin{lem}\label{inequalitysecond} For the $q_{n-1}$ case, on the whole surface $\Sigma$, if $n=2m$
\begin{eqnarray*}
&&h_1h_2|q_{n-1}|^2<h_1^{-1}h_2,\\
&&h_1h_2|q_{n-1}|^2+h_1^{-1}h_2<h_2^{-1}h_3<\dots<h_{k-1}^{-1}h_k<\dots<h_{m-1}^{-1}h_{m}<h_{m}^{-2}.
\end{eqnarray*}
If $n$ is odd, the last term is replaced by $h_{m}^{-1}.$
\end{lem}
For the $q_n$ case, the $(1,1)$-entry of $[\phi,\phi^{*_{h}}]$ is $h_1^2|q_{n}|^2-h_1^{-1}h_2$. It is strictly negative everywhere by Lemma \ref{inequalityfirst}. For the $q_{n-1}$ case, the $(1,1)$-entry of $[\phi,\phi^{*_{h}}]$ is $h_1h_2|q_{n-1}|^2-h_1^{-1}h_2$. It is strictly negative everywhere by Lemma \ref{inequalitysecond}. Therefore Part (1) follows.

For part (2), the sectional curvature in the symmetric space $G/K$ for the $2$-dimensional subspace spanned by $Y,Z$ is \begin{\eq}
 K_{G/K}(Y,Z):=-\frac{B([Y_0,Z_0],[Y_0,Z_0])}{B(Y_0,Y_0)B(Z_0,Z_0)}=-\frac{|[Y_0,Z_0]|^2}{|Y_0|^2|Z_0|^2}
\end{\eq}
where $Y_0=\Phi(Y)\in \lik^{\bot}$ and $Z_0=\Phi(Z)\in \lik^{\bot}$ from Section \ref{pre}.

The tangent space of the minimal immersion at $f(p)$ inside $G/K$ is spanned by $Y=f_*(\frac{\partial}{\partial x})$ and $Z=f_*(\frac{\partial}{\partial y})$ with $Y_0=\Phi(\frac{\partial}{\partial x})=(\phi+\phi^*)(\frac{\partial}{\partial x})=\phi(\frac{\partial}{\partial z})+\phi^*(\frac{\partial}{\partial \bar{z}})$, and $Z_0=\Phi(\frac{\partial}{\partial y})=(\phi+\phi^*)(\frac{\partial}{\partial y})=\sqrt{-1}\phi(\frac{\partial}{\partial z})-\sqrt{-1}\phi^*(\frac{\partial}{\partial \bar{z}})$. Hence \begin{\eq}
[Y_0,Z_0]=-2\sqrt{-1}[\phi(\frac{\partial}{\partial z}),\phi^*(\frac{\partial}{\partial\bar z})]=-2\sqrt{-1}[\phi,\phi^*](\frac{\partial}{\partial z}, \frac{\partial}{\partial\bar z})
\end{\eq}
By part (1), $[\phi,\phi^*]\neq 0$, hence $K_{G/K}<0.$ Therefore Part (2) follows.

Part (3) directly follows from Lemma \ref{inequalityfirst}, \ref{inequalitysecond}.
\qed\\

\textit{Proof of Lemma \ref{inequalityfirst}.}
We focus on the case $n=2m$. Along the proof, we single out each time the only differences for the case $n=2m+1$ and show that the proof works as well in this case. Let $h_k=a_kg_0^{k-\frac{n+1}{2}},  1\leq k\leq m$. Let $|q_{n}|^{2}=a_0g_0^{n}$. Then $a_{k},  0\leq k\leq m$ are globally defined functions.

The situation $q_n=0$ is clear. We then assume $q_n$ is not identically zero. Since $q_n$ is holomorphic, $q_n$ only has discrete zeros. Denote $U=\{p\in \Sigma| q_{n}(p)\neq 0\}$. Then $a_k$ satisfies, locally
\begin{\eq}
\triangle \log a_1+(1-\frac{n+1}{2}+a_1^{-1}a_2-a_{1}^{2}a_0)g_{0}&=&0,\\
\triangle \log a_k+(k-\frac{n+1}{2}+a_k^{-1}a_{k+1}-a_{k-1}^{-1}a_{k})g_0&=&0,\quad 2\leq k\leq m-1\\
\triangle \log a_m+(m-\frac{n+1}{2}+a_m^{-2}-a_{m-1}^{-1}a_{m})g_0&=&0.
\end{\eq}
 For the case $n=2m+1$, in the last equation, $a_m^{-2}$ is replaced by $a_m^{-1}$.

Let $\Omega_{k}=\log a_{k}$, $k=0,1,\dots,m$. Notice that at the zero point of $q_{n}$, $\Omega_{0}$ goes to $-\infty$ continuously. Then the system of equations becomes,
\begin{\eq}
\triangle \Omega_1+(1-\frac{n+1}{2}+e^{-\Omega_{1}+\Omega_{2}}-e^{2\Omega_{1}+\Omega_{0}})g_{0}&=&0, \quad \text{in $U$}\\
\triangle \Omega_k+(k-\frac{n+1}{2}+e^{-\Omega_{k}+\Omega_{k+1}}-e^{-\Omega_{k-1}+\Omega_{k}})g_0&=&0,\quad 2\leq k\leq m-1\\
\triangle \Omega_m+(m-\frac{n+1}{2}+e^{-2\Omega_{m}}-e^{-\Omega_{m-1}+\Omega_{m}})g_0&=&0.
\end{\eq}
 For the case $n=2m+1$, in the last equation, $e^{-2\Omega_{m}}$ is replaced by $e^{-\Omega_{m}}$.

Let $f_{1}=2\Omega_{1}+\Omega_{0}$, $f_{k}=-\Omega_{k-1}+\Omega_{k}$, $2\leq k\leq m$, $f_{m+1}=-2\Omega_{m}$. For the case $n=2m+1$, let $f_{m+1}=-\Omega_{m}$. Notice that, since $q_n$ is holomorphic, away from zeros of $q_{n}$,
\begin{\eq}
\triangle \Omega_{0}&=&\triangle \log a_{0}=\triangle\log (|q_n|^2g_{0}^{-n})\\
&=&\triangle \log|q_n|^2-n\triangle \log g_0=-ng_0
\end{\eq}
Then we obtain
\begin{\eq}
\triangle f_1+(1+2e^{f_2}-2e^{f_1})g_0&=&0, \quad \text{in $U$}\\
\triangle f_k+(1+e^{f_{k+1}}-2e^{f_{k}}+e^{f_{k-1}})g_0&=&0,\quad 2\leq k\leq m\\
\triangle f_{m+1}+(1-2e^{f_{m+1}}+2e^{f_{m}})g_0&=&0.
\end{\eq}
For the case $n=2m+1$, the last equation is replaced by \begin{\eq}
\triangle f_{m+1}+(1-e^{f_{m+1}}+e^{f_{m}})g_0&=&0.\end{\eq}
Then 
\begin{eqnarray}\label{inequality}
\triangle (f_1-f_2)+(-3(e^{f_1}-e^{f_2})+(e^{f_2}-e^{f_3}))g_0&=&0, \quad \text{in $U$}\label{contra}\\
\triangle (f_k-f_{k+1})+(-2(e^{f_k}-e^{f_{k+1}})+(e^{f_{k-1}}-e^{f_{k}})+(e^{f_{k+1}}-e^{f_{k+2}}))g_0&=&0,\quad 2\leq k\leq m-1\nonumber\\
\triangle (f_m-f_{m+1})+(-3(e^{f_m}-e^{f_{m+1}})+(e^{f_{m-1}}-e^{f_{m}}))g_0&=&0.\nonumber
\end{eqnarray}
For the case $n=2m+1$, the last equation is replaced by \begin{\eq}
\triangle (f_m-f_{m+1})+(-2(e^{f_m}-e^{f_{m+1}})+(e^{f_{m-1}}-e^{f_{m}}))g_0&=&0.\end{\eq}
Let $A_{k}$ be the maximum of $e^{f_k-f_{k+1}}$, $1\leq k\leq m$.

Step 1: We first show $A_m\leq 1$. From the compactness of $\Sigma$, $A_{k}$, $2\leq k\leq m$ must be achieved at some $p_{k}$. For $A_{1}$, since $f_{1}-f_{2}$ goes to $-\infty$ around the zeros of $q_n$, $A_1$ must be achieved at some $p_{1}\in U$. Notice that $A_k\geq 0$, $1\leq k\leq m$. Then by the maximum principle, at $p_{1}$
\begin{\eq}
&&-3(e^{f_1(p_1)}-e^{f_2(p_1)})+(e^{f_2(p_1)}-e^{f_3(p_1)})\geq 0\\
\Rightarrow&&-3(e^{(f_1-f_2)(p_1)}-1)+(1-e^{-(f_2-f_3)(p_1)})\geq 0\\
\Rightarrow&&-3(A_1-1)+(1-A_{2}^{-1})\geq 0\\
\Rightarrow&&2(A_1-1)\leq (1-A_{2}^{-1})-(A_1-1).
\end{\eq}
At $p_k$, $2\leq k\leq m-1$
\begin{\eq}
&&-2(e^{f_k(p_k)}-e^{f_{k+1}(p_k)})+(e^{f_{k-1}(p_k)}-e^{f_{k}(p_k)})+(e^{f_{k+1}(p_k)}-e^{f_{k+2}(p_k)})\geq 0\\
\Rightarrow&&-2(e^{(f_k-f_{k+1})(p_k)}-1)+(e^{(f_{k-1}-f_{k})(p_k)}-1)e^{(f_{k}-f_{k+1})(p_k)}
+(1-e^{-(f_{k+1}-f_{k+2})(p_k)})\geq 0\\
\Rightarrow&&-2(A_k-1)+(A_{k-1}-1)A_{k}+(1-A_{k+1}^{-1})\geq 0\\
\Rightarrow&&(1-A_{k}^{-1})-(A_{k-1}-1)\leq A_{k}^{-1}((1-A_{k+1}^{-1})-(A_{k}-1)).
\end{\eq}
Similarly, at $p_m$, we obtain
\begin{\eq}
(1-A_{m}^{-1})-(A_{m-1}-1)\leq -2(1-A_{m}^{-1}).
\end{\eq}
For $n=2m+1$ case, at $p_m$, instead we obtain
\begin{\eq}
(1-A_{m}^{-1})-(A_{m-1}-1)\leq -(1-A_{m}^{-1}).
\end{\eq}
Let $B_1=2(A_1-1)$, $B_{k}=(1-A_{k}^{-1})-(A_{k-1}-1)$, $2\leq k\leq m$. Then
\begin{\eq}
B_1&\leq&B_2,\\
B_{k}&\leq& A_{k}^{-1}B_{k+1}, \quad 2\leq k\leq m-1\\
B_{m}&\leq& -2(1-A_{m}^{-1}).
\end{\eq}
For $n=2m+1$ case, the last inequality is replaced by \begin{\eq}
B_{m}&\leq& -(1-A_{m}^{-1}).\end{\eq}

On the other hand,
\begin{\eq}
1-A_{m}^{-1}&=&B_{m}+(A_{m-1}-1)\\
&=&B_{m}+A_{m-1}(1-A_{m-1}^{-1})\\
&=&B_{m}+A_{m-1}(B_{m-1}+(A_{m-2}-1))\\
&\vdots&\\
&=&B_{m}+A_{m-1}B_{m-1}+A_{m-1}A_{m-2}B_{m-2}+\cdots\\
&&+A_{m-1}A_{m-2}\cdots A_{2}B_{2}+A_{m-1}A_{m-2}\cdots A_{2}(A_1-1)\\
&\leq&\underbrace{B_{m}+B_{m}+\cdots+B_{m}}_{\text{$m-1$ terms}}+\frac{1}{2}B_m\\
&\leq& -2(m-\frac{1}{2})(1-A_{m}^{-1}).
\end{\eq}
For $n=2m+1$ case, the inequality is replaced by \begin{\eq}1-A_{m}^{-1}\leq -(m-\frac{1}{2})(1-A_{m}^{-1}).\end{\eq}
Thus we obtain $A_{m}\leq 1$.

Step 2: Now we claim $A_1<1$ and $A_m<1$. To prove $A_1<1$ by contradiction, assume $A_1\geq 1$. Then $B_1\geq 0$. And then $B_k\geq 0$ for $2\leq k\leq m$. In other words, $1-A_{k}^{-1}\geq A_{k-1}-1$ for $2\leq k\leq m$. From $A_1\geq 1$, we see $A_k\geq 1$ for $2\leq k\leq m$. In particular, $A_m\geq 1$. So $A_m$ has to be equal to $1$. It turns out $A_k=1$ for every $k$. We use the strong maximum principle to get contradiction. Consider equation (\ref{contra})
\begin{\eq}
\triangle (f_1-f_2)+(-3(e^{f_1-f_2}-1)+(1-e^{-(f_2-f_3)}))e^{f_2}g_{0}&=&0, \quad \text{in $U$}
\end{\eq}
Since $e^{f_2-f_3}\leq 1$, we have
\begin{\eq}
\triangle (f_1-f_2)-3(e^{f_1-f_2}-1)e^{f_2}g_{0}\geq 0\quad \text{in $U$}.
\end{\eq}
Let $u=f_1-f_2$. Then
\begin{\eq}
\triangle u-(3e^{f_2}g_{0}\int_{0}^{1}e^{tu}dt)u\geq 0\quad \text{in $U$}.
\end{\eq}
By the assumption, $u$ achieves its maximum $0$ in $U$. Then by the strong maximum principle, $u$ has to be $0$ identically in $U$. But $q_n$ must have zeros. So $u$ has to take value $-\infty$, which is a contradiction to $u\equiv 0$. So we have proved $A_m<1$ by contradiction. But since $A_1\geq 1$ implies $A_m\geq 1$, we also see $A_1<1$ by contradiction as well, and $e^{f_1}<e^{f_2}$ everywhere on $\Sigma$.\\

Step 3: Now we begin to show $A_k<1$, for all $1\leq k\leq m$.
Assume there is some $1\leq k_0\leq m-1$ such that
$B_{k_0}\leq 0, B_{k_0+1}>0.$ The proof below works in the other two cases $B_1>0$ and $B_m\leq 0$ as well. Since $B_k\leq A_k^{-1} B_{k+1}$, for all $1\leq k\leq m-1$, then
\begin{\eq}
B_{k}&\leq&0, \quad 1\leq k\leq k_0\\
B_{k}&>&0, \quad k_0+1\leq k\leq m.
\end{\eq}
Then $B_{k}=(1-A_{k}^{-1})-(A_{k-1}-1)\leq 0~ \text{for~}1\leq k\leq k_0$ gives \begin{\eq}1-A_{k}^{-1}\leq A_{k-1}-1,  \quad 1\leq k\leq k_0. \end{\eq}
Since $A_1<1$, then $A_1-1<0$, hence \begin{\eq}A_k-1<0 , \quad 1\leq k\leq k_0.\end{\eq}
And $B_{k}=(1-A_{k}^{-1})-(A_{k-1}-1)>0~ \text{for~} k_0+1\leq k\leq m$ gives \begin{\eq}1-A_{k+1}^{-1}>A_{k}-1,\quad k_0+1\leq k\leq m.\end{\eq}
Since $A_m-1<0$, we have \begin{\eq}A_k-1<0 , \quad k_0+1\leq k\leq m.\end{\eq}
Therefore, we obtain $A_k<1,$ for all $1\leq k\leq m$.
\qed\\

\textit{Proof of Lemma \ref{inequalitysecond}.}
We only prove for the case $n=2m$. The case for $n=2m+1$ is proved in a similar way as in Lemma \ref{inequalityfirst}. As above, let $h_k=a_kg_0^{k-\frac{1}{2}-m},  1\leq k\leq m$. Let $|q_{n-1}|^{2}=a_0g_0^{n-1}$. Then $a_k$ satisfies, locally
\begin{\eq}
\triangle \log a_1+(\frac{1}{2}-m+a_1^{-1}a_2-a_{1}a_{2}a_0)g_{0}&=&0,\\
\triangle \log a_2+(\frac{3}{2}-m+a_2^{-1}a_3-a_1^{-1}a_2-a_{1}a_{2}a_0)g_{0}&=&0,\\
\triangle \log a_k+(k-\frac{1}{2}-m+a_k^{-1}a_{k+1}-a_{k-1}^{-1}a_{k})g_0&=&0,\quad 3\leq k\leq m-1\\
\triangle \log a_m+(-\frac{1}{2}+a_m^{-2}-a_{m-1}^{-1}a_{m})g_0&=&0.
\end{\eq}
Denote $U=\{p\in \Sigma| q_{n-1}(p)\neq 0\}$. Let $\Omega_{k}=\log a_{k}$, $k=0,1,\dots,m$. Then
\begin{\eq}
\triangle \Omega_1+(\frac{1}{2}-m+e^{-\Omega_{1}+\Omega_{2}}-e^{\Omega_{1}+\Omega_2+\Omega_{0}})g_{0}&=&0, \quad \text{in $U$}\\
\triangle \Omega_2+(\frac{3}{2}-m+e^{-\Omega_{2}+\Omega_{3}}-e^{-\Omega_{1}+\Omega_{2}}-e^{\Omega_{1}+\Omega_{2}+\Omega_{0}})g_{0}&=&0, \quad \text{in $U$}\\
\triangle \Omega_k+(k-\frac{1}{2}-m+e^{-\Omega_{k}+\Omega_{k+1}}-e^{-\Omega_{k-1}+\Omega_{k}})g_0&=&0,\quad 3\leq k\leq m-1\\
\triangle \Omega_m+(-\frac{1}{2}+e^{-2\Omega_{m}}-e^{-\Omega_{m-1}+\Omega_{m}})g_0&=&0.
\end{\eq}
Let $f_{1}=\Omega_{1}+\Omega_{2}+\Omega_{0}$, $f_{k}=-\Omega_{k-1}+\Omega_{k}$, $2\leq k\leq m$, $f_{m+1}=-2\Omega_{m}$. Notice that $\triangle \Omega_{0}=(1-n)g_0$. Then we obtain
\begin{\eq}
\triangle f_1+(1+e^{f_3}-2e^{f_1})g_{0}&=&0, \quad \text{in $U$}\\
\triangle f_2+(1+e^{f_3}-2e^{f_2})g_{0}&=&0,\\
\triangle f_3+(1+e^{f_4}-2e^{f_3}+e^{f_2}+e^{f_1})g_{0}&=&0,\\
\triangle f_k+(1+e^{f_{k+1}}-2e^{f_{k}}+e^{f_{k-1}})g_0&=&0,\quad 4\leq k\leq m\\
\triangle f_{m+1}+(1-2e^{f_{m+1}}+2e^{f_{m}})g_0&=&0.
\end{\eq}

Then we have the following system
\begin{eqnarray*}
\triangle (f_3-f_{4})+(-2(e^{f_3}-e^{f_{4}})+(e^{f_{4}}-e^{f_5})+(e^{f_1}+e^{f_2}-e^{f_{3}}))g_0&=&0,\\
\triangle (f_k-f_{k+1})+(-2(e^{f_k}-e^{f_{k+1}})+(e^{f_{k-1}}-e^{f_{k}})+(e^{f_{k+1}}-e^{f_{k+2}}))g_0&=&0,\quad 3\leq k\leq m-1\nonumber\\
\triangle (f_m-f_{m+1})+(-3(e^{f_m}-e^{f_{m+1}})+(e^{f_{m-1}}-e^{f_{m}}))g_0&=&0.\nonumber
\end{eqnarray*}
Let $A_{k}$ be the maximum of $e^{f_k-f_{k+1}}, k=1$ and $3\leq k\leq m$, and $A_2$ be the maximum of $e^{f_1-f_3}+e^{f_2-f_3}.$

Step 1: We show $A_1<1.$
We only look at the first equation of the Hitchin equation
\begin{\eq}
\triangle \Omega_1+(\frac{1}{2}-m+e^{-\Omega_{1}+\Omega_{2}}-e^{\Omega_{1}+\Omega_2+\Omega_{0}})g_{0}&=&0, \quad \text{in $U$}
\end{\eq}
Consider the function $-\frac{1}{2}\Omega_0=\log|q_{n-1}|^{-1}g_0^{\frac{n-1}{2}}$, it satisfies the above equation. Since it is $\infty$ at zeros of $q_{n-1}$, it gives an upper bound for $\Omega_1$. By the strong maximum principle, we have $\Omega_1<\log|q_{n-1}|^{-1}g_0^{\frac{n-1}{2}}.$
Hence $e^{2\Omega_1}|q_{n-1}|^2<1$, then $A_1<1.$


Step 2: We claim that $A_m<1$. For other $A_k$, we apply the same argument in Lemma \ref{inequalityfirst} to the above system.
 Let $B_{k}=(1-A_{k}^{-1})-(A_{k-1}-1)$, $3\leq k\leq m$. Then
\begin{eqnarray}\label{Bk}
B_{k}&\leq& A_{k}^{-1}B_{k+1}, \quad 3\leq k\leq m-1\\
B_{m}&\leq& -2(1-A_{m}^{-1}).\nonumber
\end{eqnarray}

Let $B_2=A_2-1$. Now we start to obtain inequality for $B_2$.
Applying the inequality
\begin{\eq}
\triangle\log(e^f+1)\geq \frac{e^f}{e^f+1}\triangle f,
\end{\eq}
we have \begin{\eq}
\frac{1}{g_0}\triangle\log(e^{f_1-f_3}+e^{f_2-f_3})&=&\triangle\log(e^{f_1-f_2}+1)+\triangle\log e^{f_2-f_3}\\
&\geq&\frac{e^{f_1-f_2}}{e^{f_1-f_2}+1}\triangle(f_1-f_2)+\triangle(f_2-f_3)\\
&=&\frac{e^{f_1-f_2}}{e^{f_1-f_2}+1}2(e^{f_1}-e^{f_2})+(e^{f_1}+3e^{f_2}-3e^{f_3}+e^{f_4})\\
&=& e^{f_2}\frac{(e^{f_1-f_2}-1)^2}{e^{f_1-f_2}+1}+(e^{f_1}-e^{f_2})+(e^{f_1}+3e^{f_2}-3e^{f_3}+e^{f_4})\\
&>&2(e^{f_1}+e^{f_2}-e^{f_3})-(e^{f_3}-e^{f_4}), \quad\text{By Step 1}\\
&=&2e^{f_3}(e^{f_1-f_3}+e^{f_2-f_3}-1)-e^{f_3}(1-e^{-f_3+f_4})
\end{\eq}
Hence, at maximum of $e^{f_1-f_3}+e^{f_2-f_3},$ we have
\begin{\eq}
2(A_2-1)<1-A_3^{-1}.
\end{\eq}
Hence \begin{equation}\label{B2}
B_2=A_2-1<(1-A_3^{-1})-(A_2-1)=B_3.\end{equation}

We apply the same argument:
\begin{\eq}
1-A_{m}^{-1}&=&B_{m}+(A_{m-1}-1)\\
&=&B_{m}+A_{m-1}(B_{m-1}+(A_{m-2}-1))\\
&&\cdots\\
&=&B_{m}+A_{m-1}B_{m-1}+A_{m-1}A_{m-2}B_{m-2}+\cdots\\
&&+A_{m-1}A_{m-2}\cdots A_{3}B_{3}+A_{m-1}A_{m-2}\cdots A_{3}(A_2-1)\\
\text{Applying Inequality (\ref{Bk}),(\ref{B2})   }&<&\underbrace{B_{m}+B_{m}+\cdots+B_{m}}_{\text{$m-1$ terms}}\\
&<& -2(m-1)(1-A_{m}^{-1}).
\end{\eq}
Therefore $A_m<1.$

Step 3: By applying the same process as in Lemma \ref{inequalityfirst}, we show $A_k<1$ and finish the proof.\qed

\section{Comparison to decoupled equations}\label{decouple}
In this section, we compare the solution to coupled Hitchin equations by the solutions to decoupled vortex-like equations. For $n$ even, set $n=2m$; for $n$ odd, set $n=2m+1$. Throughout this section, we assume $m\geq 2$, since for $m=1$ the Hitchin system of equations is a single equation which is already decoupled.
\begin{thm}\label{scalar}
There exists a unique Hermitian metric $u_k, v_k$ respectively on the holomorphic line bundle $K^{\frac{n+1-2k}{2}}$ satisfying, locally
\begin{eqnarray}
\triangle \log u_k+u_k^{-\frac{2}{n+1-2k}}-(u_k^2|q_{n}|^2)^{\frac{1}{2k-1}}&=&0,\quad 1\leq k\leq m, \label{uk}\\
\triangle \log v_k+v_k^{-\frac{2}{n+1-2k}}-(v_k^2|2q_{n-1}|^2)^{\frac{1}{2k-2}}&=&0,\quad 2\leq k\leq m. \label{vk}
\end{eqnarray}
And the following estimates hold,
\begin{eqnarray*}
\max\{|q_n|^{\frac{2\alpha_k}{n}},(\alpha_kg_0)^{\alpha_k}\}\leq &u_k^{-1}&<((\underset{\Sigma}{\max}|q_n|_{g_0}^{\frac{2}{n}}+\alpha_k)g_0)^{\alpha_k},\quad 1\leq k\leq m,\\
\max\{|2q_{n-1}|^{\frac{2\alpha_k}{n-1}}, (\alpha_kg_0)^{\alpha_k}\}\leq &v_k^{-1}&<((\underset{\Sigma}{\max}|2q_{n-1}|_{g_0}^{\frac{2}{n-1}}+\alpha_k)g_0)^{\alpha_k}, \quad 2\leq k\leq m,
\end{eqnarray*}
where $\alpha_k=\frac{n+1-2k}{2}$ and $g_0$ is the unique Hermitian hyperbolic metric on $K^{-1}$.
\end{thm}
\begin{rem}
Such equations are natural generalization of familiar vortex-like equations. In particular, in the case $n=2,k=1$, this is the harmonic equation from surface to surface with given Hopf differential $q_2$. In the case $n=3, k=1$, this gives Wang's equation \cite{Wang} and the solution is the Blaschke metric for the hyperbolic affine sphere. Dumas and Wolf \cite{DumasWolf} solved such equations for the case when $n$ is general and $k=1$ on the complex plane.
\end{rem}
\begin{rem}
Let $u_k=b_kg^{\frac{2k-n-1}{2}}$, then the equation (\ref{uk}) becomes
\begin{\eq}
\triangle_g \log b_k-{\frac{2k-n-1}{2}} K_g+b_k^{-\frac{2}{2(m-k)+1}}-(b_k^2|q_n|_g^2)^{\frac{1}{2k-1}}=0,
\end{\eq}
where $K_{g}=-\triangle_{g} \log g$ is the sectional curvature of $g$.
In particular, there is a geometric interpretation for the equation (\ref{uk}). Let $\sigma_k=u_k^{-\frac{2}{n+1-2k}}$. Then $\sigma_k$ is an Hermitian metric on the Riemann surface $\Sigma$. The equation for $\sigma_k$ is
\begin{\eq} \frac{n+1-2k}{2}K_{\sigma_k}=-1+|q_{n}|_{\sigma_k}^{\frac{2}{2k-1}}.
\end{\eq}
\end{rem}
\textit{Proof of Theorem \ref{scalar}.} Again, we only prove the case $n=2m$ for the $q_n$ case. The proof for other cases is similar.
Take the Hermitian hyperbolic metric $g_0$ on $\Sigma$. Let $\sigma_k=u_k^{-\frac{2}{n+1-2k}}$ and $\sigma_k=g_0e^{\eta_k}$ with $\eta_k$ being a $C^{\infty}$ function on $\Sigma$.
Then the equation (\ref{uk}) is equivalent to
\begin{\eq}
\frac{n+1-2k}{2}\triangle_{g_0} \eta_k-e^{\eta}+|q_{n}|_{g_0}^{\frac{2}{2k-1}}e^{-\frac{n+1-2k}{2k-1}\eta_k}+\frac{n+1-2k}{2}=0,\quad 1\leq k\leq m.
\end{\eq}
Then showing the first part of Theorem \ref{scalar} reduces to show the existence and uniqueness of $\eta_k$ which follows from the following lemma. The proof is similar to Proposition 4.0.2 in Loftin \cite{LoftinAffine}.
\begin{lem} \label{lof}
Let $(M,g)$ be a closed Riemannian manifold, $f$ be a nonnegative $C^{\infty}$ function on $M$, and $a,b,c$ be positive real numbers. Then the equation
\begin{\eq}a\triangle_{g} \eta-e^{\eta}+f(x)e^{-c\eta}+b=0,\quad 1\leq k\leq m,\end{\eq}
has a unique $C^{\infty}$ solution, such that
\begin{\eq}
\log b\leq\eta<\log (G^{\frac{1}{c+1}}+b),
\end{\eq}
where $G$ denotes the maximum value of $f(x)$. If $f$ is not identically zero, then $\eta>\log b$.
\end{lem}
\textit{Proof of Lemma \ref{lof}.} For the existence, it is sufficient to find a subsolution and a supersolution for this equation (see Schoen and Yau \cite{SchoenYau}). First, the constant function $\eta=\log b$ is a subsolution by direct calculation. Second, set $m$ to be the smallest positive root of the equation
\begin{\eq}
x^{c+1}-bx^c-G=0.
\end{\eq}
Then $\eta=\log m$ satisfies
\begin{\eq}
a\triangle_{g} \eta-e^\eta+f(x)e^{-c\eta}+b=-m+f(x)m^{-c}+b\leq 0.
\end{\eq}
Therefore $s=\log m$ is a supersolution. Then there is a smooth solution $\eta$  to the equation satisfying
\begin{\eq}
\log b\leq\eta\leq \log m<\log (G^{\frac{1}{c+1}}+b).
\end{\eq}
If $f$ is not zero, the strong maximum principle implies $\eta>\log b$.
The uniqueness comes from the maximum principle.
And by the standard elliptic theory, we obtain the smoothness of the solution.
\qed

Theorem \ref{scalar} then follows from the above Lemma directly except to show $u_k\leq |q_n|^{-\frac{n+1-2k}{n}}$. In fact, we observe that $|q_n|^{-\frac{n+1-2k}{n}}$ is also a solution to the equation (\ref{uk}) outside the zero locus of $q_n$. Then applying the strong maximum principle, we finish the proof.
\qed\\

We now compare the solution to the coupled Hitchin system of equations to such Hermitian metrics $u_k$ and $v_k$ satisfying decoupled equations.
\begin{thm}\label{comparison}
For $m\geq 2$,
\begin{\eq}
\text{$q_{n}$ case:}& \quad h_k<u_k, \quad 1\leq k\leq m,\\
\text{$q_{n-1}$ case:}& \quad h_k< v_k, \quad 2\leq k\leq m.
\end{\eq}
\end{thm}

\begin{rem}
If $m=1$, for the $q_n$ case, then the Hitchin equation coincides with the equation for $u_1$. Therefore $h_1=u_1$. In this sense, $u_k$ is a natural generalization to bound the coupled Hitchin system.
\end{rem}
\textit{Proof of Theorem \ref{comparison}.} Again, we only prove the case $n=2m$ for the $q_n$ case. The proof for the other cases is similar. From Lemma \ref{inequalityfirst}
\begin{\eq}
(h^{-1}_k h_{k+1})^{2(m-k)+1}<(h^{-1}_k h_{k+1})^2\cdot(h_{k+1}^{-1}h_{k+2})^2\cdots(h_{m-1}^{-1}h_{m})^2\cdot h_m^{-2}=h_k^{-2},
\end{\eq}
hence
\begin{\eq}
h^{-1}_k h_{k+1}<h_k^{-\frac{2}{2(m-k)+1}},  \quad 1\leq k\leq m-1.
\end{\eq}
From Lemma \ref{inequalityfirst}\begin{\eq}
(h_{k-1}^{-1}h_{k})^{2k-1}>(h_1^2|q_{2m}|^2)(h_1^{-1}h_{2})^2\cdots(h_{k-1}^{-1}h_{k})^2=h_k^{2}|q_{2m}|^2,\end{\eq}
hence
\begin{\eq}
h_{k-1}^{-1}h_{k}>({h_k}^2|q_{2m}|^2)^{\frac{1}{2k-1}},\quad 2\leq k\leq m.
\end{\eq}
Therefore, the system of equations \begin{eqnarray*}
\triangle \log h_1+h_1^{-1}h_2-h_1^2|q_{n}|^2&=&0,\nonumber\\
\triangle \log h_k+h_k^{-1}h_{k+1}-h_{k-1}^{-1}h_{k}&=&0,\quad 2\leq k\leq m-1\\
\triangle \log h_m+h_m^{-2}-h_{m-1}^{-1}h_{m}&=&0,\nonumber
\end{eqnarray*} implies that
\begin{\eq}
\triangle \log h_k+h_k^{-\frac{2}{2(m-k)+1}}-(h_k^2|q_{2m}|^2)^{\frac{1}{2k-1}}>0,\quad 1\leq k\leq m
\end{\eq}
Then applying the strong maximum principle, we obtain $h_k<u_k$.
\qed
\begin{rem}
For $h_1$ in the $q_{n-1}$ case, from Lemma \ref{inequalitysecond}, one can show that $h_1<|q_{n-1}|^{-1}$.
\end{rem}
Now fixing $q_n$ (or $q_{n-1}$), we study the asymptotic behavior of $h_k,u_k,v_k$ along the ray $tq_n$ (or $tq_{n-1}$) when $t$ approaches to infinity.
For $h_k$, Collier and Li \cite{CL14} proved the following result.

\begin{thm}\label{CollierLi}(Collier-Li \cite{CL14}) Suppose $q_n$ (or $q_{n-1}$) is not zero. Then for each compact set $K\subset \Sigma$ away from zeros of $q_n$ (or $q_{n-1}$), there is a positive constant $C=C(K)$ independent of $t$ such that at each point $p\in K$,
\begin{\eq}
\text{$q_{n}$ case:}&&\quad h_k=|q_n|^{-\frac{n+1-2k}{n}}(1+O(||q_n||^{-1})), \quad 1\leq k\leq m,\\
\text{$q_{n-1}$ case:}&&\quad h_k=(|2q_{n-1}|)^{-\frac{n+1-2k}{n-1}}(1+O(||q_{n-1}||^{-1})),\quad 2\leq k\leq m,
\end{\eq}
where $||q_n||=\int_\Sigma|q_n|^{\frac{2}{n}}$ and $f\in O(||q_n||^{-1})$ means $|f|/||q_n||^{-1}\leq C$.
\end{thm}
\begin{rem}
By Theorem \ref{scalar} and Theorem \ref{comparison}, $h_k<|q_n|^{-\frac{n+1-2k}{n}}$, so we improve the above result in one direction.
\end{rem}
Let $u_k^t$ (or $v_k^t$) be the solution for $tq_n$ (or $tq_{n-1}$) in Theorem \ref{scalar} (or $tq_{n-1}$). Combining Theorem \ref{comparison} and Theorem \ref{CollierLi}, we obtain
\begin{thm}
 Suppose $q_n$ (or $q_{n-1}$) is not zero. Then for each compact set $K\subset \Sigma$ away from zeros of $q_n$ (or $q_{n-1}$), there is a positive constant $C=C(K)$ independent of $t$ such that
\begin{\eq}
\text{$q_{n}$ case:}&& u_k^t(1-C||tq_n||^{-1})\leq h_k^t<u_k^t, \quad 1\leq k\leq m,\\
\text{$q_{n-1}$ case:}&& v_k^t(1-C||tq_{n-1}||^{-1})\leq h_k^t<v_k^t, \quad 2\leq k\leq m.
\end{\eq}
\end{thm}

\begin{rem}
The theorem above shows that $u_k$ (or $v_k$) is an upper approximation of $h_k$. 
\end{rem}

Finally, we prove the following asymptotic estimates of $u_k,v_k$. It can be derived from Theorem \ref{scalar}, Theorem \ref{comparison} and Theorem \ref{CollierLi}. Here we give a direct proof not relying on Theorem \ref{CollierLi}. 
\begin{thm}\label{asym}
Suppose $q_n$ (or $q_{n-1}$) is not zero. Then for each compact set $K\subset \Sigma$ away from zeros of $q_n$ (or $q_{n-1}$), there is a positive constant $C=C(K)$ independent of $t$ such that
\begin{\eq}
|tq_n|^{-\frac{n+1-2k}{n}}(1-C||tq_n||^{-1})\leq &u_k^t&<|tq_n|^{-\frac{n+1-2k}{n}}, \quad 1\leq k\leq m,\\
|2tq_{n-1}|^{-\frac{n+1-2k}{n-1}}(1-C||tq_{n-1}||^{-1})\leq &v_k^t&<|2tq_{n-1}|^{-\frac{n+1-2k}{n-1}}, \quad 2\leq k\leq m.
\end{\eq}
\end{thm}

\textit{Proof of Theorem \ref{asym}.}
For simplicity, we drop the superscript $t$ of $u_k^t,v_k^t$. We only prove $u_k\geq |tq_n|^{-\frac{n+1-2k}{n}}(1-C||tq_n||^{-1})$. Let $K$ be a compact subset of $\Sigma$ which does not contain any zeros of $q_n$. Choose a background metric $g$ on the surface $\Sigma$ defined as follows: $g=\frac{|q_n|^{\frac{2}{n}}}{||q_n||}$ on $K$ and $g\geq \frac{|q_n|^{\frac{2}{n}}}{||q_n||}$ outside $K$. Let $u_k=b_kg^{\frac{2k-n-1}{2}}$, where $b_k$ is a positive function on $\Sigma$ satisfying
\begin{\eq}
\triangle_g \log b_k-{\frac{2k-n-1}{2}} K_g+b_k^{-\frac{2}{n+1-2k}}-(b_k^2|tq_n|_g^2)^{\frac{1}{2k-1}}=0.
\end{\eq}

At the minimum of $b_k$,
\begin{\eq}
\frac{n+1-2k}{2}K_g+b_k^{-\frac{2}{n+1-2k}}-(b_k^2|tq_n|_g^2)^{\frac{1}{2k-1}}\leq 0.
\end{\eq}
Let $M=\underset{\Sigma}{\max}\frac{n+1-2k}{2}|K_g|$ and hence
\begin{\eq}
-M+b_k^{-\frac{2}{n+1-2k}}\leq(b_k^2|tq_n|_g^2)^{\frac{1}{2k-1}}\leq b_k^{\frac{2}{2k-1}}||tq_n||^{\frac{n}{2k-1}}.
\end{\eq}
Let $x=||tq_n|| b_k^{\frac{2}{n+1-2k}}$, then $x$ satisfies
\begin{\eq}
x^{\frac{n}{2k-1}}+\frac{M}{||tq_n||}x-1\geq 0.
\end{\eq}
Therefore $||tq_n|| b_k^{\frac{2}{n+1-2k}}=x\geq 1-C||tq_n||^{-1}$ at the minimum of $b_k$, also the minimum of $x$.

Hence globally on the surface, $||tq_n|| b_k^{\frac{2}{n+1-2k}}=x\geq 1-C||tq_n||^{-1}$. Then
\begin{eqnarray*}
u_k&=&(x||tq_n||^{-1})^{\frac{n+1-2k}{2}}g^{-\frac{n+1-2k}{2}}\\
&\geq&x^{\frac{n+1-2k}{2}}|tq_n|^{-\frac{n+1-2k}{n}}\\
&\geq&|tq_n|^{-\frac{n+1-2k}{n}}(1-C||tq_n||^{-1}).
\end{eqnarray*}
\qed
\begin{rem}
It is an interesting question that how $C$ varies when $K$ approaches to $\Sigma-Z$, where $Z$ is the set of zeros of $q_n$.
\end{rem}

\end{document}